\documentclass[10pt]{amsart}
\usepackage{latexsym, amsmath}

\newcommand{\cal}{\mathcal}

\renewcommand{\epsilon}{\varepsilon}
\newcommand{\newsection}[1]
{\subsection{#1}\setcounter{theorem}{0} \setcounter{equation}{0}
\par\noindent}

\newtheorem{theorem}{Theorem}

\newtheorem{lemma}[theorem]{Lemma}
\newtheorem{corr}[theorem]{Corollary}

\newtheorem{proposition}[theorem]{Proposition}
\newtheorem{deff}[theorem]{Definition}

\newcommand{\bth}{\begin{theorem}}
\newcommand{\ble}{\begin{lemma}}
\newcommand{\bcor}{\begin{corr}}
\newcommand{\ltrt}{{L^2({\mathbb R}^3)}}
\newcommand{\ltfat}{{L^2({R/4<|x|<2R})}}
\newcommand{\holt}{{\dot{H}^1\times L^2({\mathbb R}^3)}}
\newcommand{\ltoo}{{L^2({\mathbb{R}}^3\backslash\mathcal{K})}}
\newcommand{\bdeff}{\begin{deff}}
\newcommand{\lirt}{{L^\infty({\mathbb R}^3)}}
\newcommand{\lsrt}{{L^6({\mathbb R}^3)}}
\newcommand{\lioo}{{L^{\infty}({\mathbb{R}}^3 \backslash \mathcal{K})}}
\newcommand{\bprop}{\begin{proposition}}
\newcommand{\eth}{\end{theorem}}
\newcommand{\ele}{\end{lemma}}
\newcommand{\ecor}{\end{corr}}
\newcommand{\edeff}{\end{deff}}

\newcommand{\eprop}{\end{proposition}}

\newcommand{\cd}{\, \cdot\, }

\newcommand{\R}{{\mathbb R}}

\newcommand{\e}{\varepsilon}

\newcommand{\supp}{\text{supp }}
\renewcommand{\Pi}{\varPi}

\renewcommand{\epsilon}{\varepsilon}

\begin{document}

\title[Almost global existence for quasilinear wave equations]
{Almost global existence for  quasilinear wave
equations in three space dimensions}
\thanks{The authors were supported in part by the NSF}
\thanks{
The authors would like to thank the referee for a number of
helpful suggestions and comments.  In particular, we are grateful
for the suggestion that non-diagonal systems should be handled.  The third author would like to thank J.-M. Delort for helpful
conversations regarding this paper.  He would also like to thank
Paris Nord for their hospitality during a recent visit where part
of this research was carried out}
\author{Markus Keel}
\author{Hart F. Smith}
\author{Christopher D. Sogge}
\address{Department of Mathematics, University of Minnesota,
Minneapolis, MN 55455}
\address{Department of Mathematics, University of Washington,
Seattle, WA 98195}
\address{Department of Mathematics, The Johns Hopkins University,
Baltimore, MD 21218}

\maketitle

\newsection{Introduction}

This article studies almost global existence for solutions of
quadratically quasilinear systems of wave equations in three space
dimensions.  The approach here uses only the classical invariance
of the wave operator under translations, spatial rotations, and
scaling.  Using these techniques we can handle wave equations in
Minkowski space or Dirichlet-wave equations in the exterior of a
smooth, star-shaped obstacle.  We can also apply our methods to
systems of quasilinear wave equations having different wave
speeds.

This extends our work \cite{KSS2} for the semilinear case.
Previous almost global existence theorems  for quasilinear
equations in three space dimensions were for the non-obstacle
case. In \cite{JK}, John and Klainerman proved almost global
existence on Minkowski space for quadratic, quasilinear equations
using the Lorentz invariance of the wave operator in addition to
the symmetries listed above. Subsequently, in \cite{KS}, Klainerman
and Sideris obtained the same result for a class of quadratic, divergence-form
nonlinearities without relying on Lorentz invariance.  This line
of thought was refined and applied to prove global in time results
for null-form equations related to the theory of elasticity
in Sideris \cite{Si, Si2}, and for multiple speed systems of
null-form quasilinear equations in  Sideris and Tu \cite{Si3}, and
Yokoyama \cite{Y}.

The main difference between our approach and the earlier ones is
that we exploit the $O(|x|^{-1})$ decay of solutions of
wave equations with sufficiently decaying initial data  as much as we involve
the stronger $O(t^{-1})$ decay.  Here, of course, $x=(x_1,x_2,x_3)$ is
the spatial component, and $t$ the time component, of a space-time
vector $(t,x)\in \R_+\times \R^3$.   Establishing
$O(|x|^{-1})$ decay is considerably easier and can be achieved
using only the invariance with respect to translations and spatial
rotation.
A weighted $L^2$ space-time estimate for inhomogeneous wave equations
(Proposition \ref{prop3.1} below, from \cite{KSS2}) is important in making
the spatial decay useful for the long-time existence argument.

For semilinear systems, one can show almost global existence from
small data using only this spatial decay \cite{KSS2}. For
quasilinear systems, however, we also have to show that both first
and second derivatives of $u$ decay like $1/t$. Fortunately, we
can do this using a variant of some $L^1\to L^\infty$ estimates of
John, H\"ormander, and Klainerman (see \cite{H}, Lemma 6.6.8, and
also \cite{h2}, \cite{knull}) that is well adapted to our approach
since it only uses the Euclidean rotation and scaling vector
fields and involves $1/|x|$ decay.

The translation, rotation, and scaling vector fields
are useful for obstacle problems since
their normal components to the boundary of the obstacle in space-time
are $O(1)$. The Lorentz boost fields, which were also
used in the original generalized energy
approach \cite{JK}, do not have this property for any obstacle: these fields
$t\partial_i+x_i\partial_t$,
$i=1,2,3$, have normal components of size $t$. Consequently it
seems difficult to use these Lorentz boosts and still obtain optimal results.

In the Minkowski space (single-speed) setting all of the
generators of the Lorentz group can be used without difficulty
just by using the fact that they have favorable commutation
properties with the D'Alembertian. In the case of an obstacle
problem, however, not even the Euclidean rotation or scaling
vector fields commute with the Dirichlet-wave operator. Because of
the the boundary conditions, the generalized energy estimates here
are more involved than they are for the Minkowski space setting,
particularly when these estimates involve the scaling vector field
$t\partial_t +x\cdot \nabla_x$.
For the scaling field we have
to use our assumption that the obstacle is star-shaped in
an argument that is reminiscent of that of Morawetz \cite{M}.

We now describe more precisely the initial boundary value problems we shall
consider.  We assume that the obstacle ${\cal K} \subset \R^3$ is smooth
and strictly star-shaped with respect to the origin.
By this, we understand that in
polar coordinates $x=r\omega$, $(r,\omega)\in [0,\infty)\times S^2$,
we can write
\begin{equation}\label{i.2}
{\cal K}=\{\,(r,\omega) \, : \, \phi(\omega)-r \ge 0\},
\end{equation}
where $\phi$ is a smooth {\it positive} function on $S^2$.  Thus,
$$0\in {\cal K}, \, \, \, \text{but } \, 0\notin \partial{\cal K}
=\{x: \, r=\phi(\omega)\}.$$ For such  ${\mathcal K}\subset \R^3$, we
consider  smooth, quadratic, quasilinear
systems of the form
\begin{equation}\label{0.2}
\begin{cases}
\square_c u = Q(du,d^2u), \quad (t,x)\in \R_+\times
\R^3\backslash {\mathcal K}
\\
u(t,\cd)|_{\mathcal K}=0
\\
u(0,\cd)=f, \, \, \partial_tu(0,\cd)=g.
\end{cases}
\end{equation}
Here
\begin{equation}\label{squarem}
\square_c = (\square_{c_1},\square_{c_2},\dots,\square_{c_N})
\end{equation} is a vector-valued multiple speed D'Alembertian
with
$$\square_{c_I}=\partial^2_t-c^2_I\Delta,$$
where we assume that the wave speeds $c_I$ are all positive but
not necessarily distinct.  Here $\Delta =
\partial_1^2+\partial_2^2+\partial_3^2$ is the standard
Laplacian.

By quasilinear we mean that the nonlinear term $Q(du,d^2u)$ is
linear in the second derivatives of $u$.  We shall also assume
that the highest order nonlinear terms are symmetric, by which we
mean that, if we let $\partial_0=\partial_t$, then
\begin{equation}\label{0.3}
Q^I(du,d^2u)=B^I(du)+\sum_{\substack{0\le j,k,l\le 3\\ 1\le J,K\le N} }
B^{IJ,jk}_{K,l}\partial_lu^K \,
\partial_j\partial_ku^J, \quad 1\le I\le N,
\end{equation}
with $B^I(du)$ a quadratic form in the gradient of $u$, and
$B^{IJ,jk}_{K,l}$ real constants satisfying the symmetry
conditions
\begin{equation}\label{symm0}
B^{IJ,jk}_{Kl}=B^{JI,jk}_{Kl}=B^{IJ,kj}_{Kl}.
\end{equation}
The second equation here places no restriction on our systems as
we may obviously ensure this by symmetrizing.  The first equality
in \eqref{symm0} will be used when we prove the standard energy
estimates.  Some restriction along these lines seems necessary for
our theorem to be true. In fact, there are even simple examples of
{\em linear} second order systems which violate  \eqref{symm0} and
for which the basic energy estimate fails.  (This failure is
well-known, for example it is pointed out by Fritz John in his
work on elasticity.) For completeness, we will sketch one such
example following Proposition \ref{prop3.2} below.

In order to solve \eqref{0.2} we must also assume that the data
satisfies the relevant compatibility conditions.  Since these are
well known (see e.g., \cite{KSS}), we shall describe them briefly.
To do so we first let $J_ku =\{\partial^\alpha_xu: \, 0\le
|\alpha|\le k\}$ denote the collection of all spatial derivatives
of $u$ of order up to $k$.  Then if $m$ is fixed and if $u$ is a
formal $H^m$ solution of \eqref{0.2} we can write
$\partial_t^ku(0,\cd)=\psi_k(J_kf,J_{k-1}g)$, $0\le k\le m$, for
certain compatibility functions $\psi_k$ which depend on the
nonlinear term $Q$ as well as $J_kf$ and $J_{k-1}g$.  Having done
this, the compatibility condition for \eqref{0.2} with $(f,g)\in
H^m\times H^{m-1}$ is just the requirement that the $\psi_k$
vanish on $\partial{\mathcal K}$ when $0\le k\le m-1$.
Additionally, we shall say that $(f,g)\in C^\infty$ satisfy the
compatibility conditions to infinite order if this condition holds
for all $m$.

We can now state our main result.  In describing the initial data we shall
use the weight
\begin{align}
\langle x \rangle & \equiv (1 + |x|^2)^{\frac{1}{2}}\,.
\end{align}

\begin{theorem}\label{theorem0.1}
Let ${\mathcal K}$ be a star-shaped obstacle, and assume that
$Q(du,d^2u)$ and $\square_c$ are as above. Assume further that
$(f,g)\in C^\infty(\R^3\backslash \mathcal{K})$ satisfies
the compatibility conditions to infinite order.

Then there are constants $\kappa, \varepsilon_0>0$, and an integer 
$N > 0$ so that for all $\varepsilon \leq \varepsilon_0$, if
\begin{equation}\label{0.7}
\sum_{|\alpha|\le N}\|\langle x\rangle^{|\alpha|} \partial_x^\alpha
f\|_\ltoo + \sum_{|\alpha|\le N-1}\|\langle
x\rangle^{|\alpha| + 1} \partial_x^\alpha g\|_\ltoo \le
\varepsilon,\end{equation} then \eqref{0.2} has a unique solution
$u\in C^\infty([0,T_\varepsilon]\times \R^3\backslash
\mathcal{K})$, with
\begin{equation}\label{0.8}
T_\varepsilon = \exp(\kappa/\varepsilon).
\end{equation}
\end{theorem}

The norms in which we control the solution up to time $T_\varepsilon$ are
found in \S \ref{section:almostglobaloutside}.
%\eqref{newMk} below.

We shall actually establish existence of limited regularity almost
global solutions $u$ for data $(f,g)\in H^{N}\times H^{N-1}$
satisfying the relevant compatibility conditions.  The fact
that $u$ must be smooth if $f$ and $g$ are smooth and satisfy the
compatibility conditions of infinite order follows from standard
local existence theorems (see \S 9, \cite{KSS}).  Also, we are not
concerned here with minimal regularity issues.  The
value $N=15$ which we eventually require (see \eqref{10.1} below) is certainly not optimal.

Together with the finite propagation speed of our equations, the
blow-up examples in e.g. John \cite{johnblowup} show that for the
class of nonlinearities described above, the  time of existence
\eqref{0.8} is sharp. If we restrict our attention to null-form
nonlinearities and single speed systems, global in time solutions
outside of star-shaped obstacles were established by the authors
in \cite{KSS}.This extended  earlier spherically 
symmetric work of \cite{godin}.
For related work outside of obstacles in higher dimensions
see  \cite{hayashi}, \cite{tsutsumi}.

We point out that results similar to those in Theorem
\ref{theorem0.1} were announced in Datti \cite{D}, but there
appears to be a gap in the argument which has not been repaired.
Specifically, the proof of Theorem 5.3 of \cite{D} cannot be attained as
claimed, and hence the main estimates of the paper remain unproven.

As we remarked before, we can also give a  proof of a multiple
speed generalization of the almost global existence theorem of
John and Klainerman \cite{JK}:

\begin{theorem}\label{theorem0.2}
Assume that $Q(du,d^2u)$  and $\square_c$ are as above. Then there
exists $N > 0$ and constants $\kappa, \varepsilon_0>0$ so that for
all $\varepsilon< \varepsilon_0$ and data $(f,g)\in
C^\infty(\R^3) \cap L^6(\R^3)$ satisfying
\begin{equation}\label{0.9}
\sum_{|\alpha|\le N}
\|\langle x\rangle^{|\alpha|} \partial_x^\alpha f'\|_\ltrt + 
\sum_{|\alpha|\le N-1}
\|\langle x \rangle^{|\alpha|}\partial_x^\alpha g\|_\ltrt \leq \varepsilon,
\end{equation}
the system,
\begin{equation}\label{0.10}
\begin{cases}
\square_c u = Q(du,d^2u)
\\
u(0,\cd)=f, \, \, \partial_tu(0,\cd)=g
\end{cases}
\end{equation}
has a unique solution $u\in C^\infty([0,T_\varepsilon]\times
\R^3)$, where
\begin{equation}\label{0.11}
T_\varepsilon = \exp(\kappa/\varepsilon).
\end{equation}
\end{theorem}

As we noted earlier, in \cite{KS} Klainerman and Sideris established Theorem
\ref{theorem0.2} in the case of certain divergence-form
nonlinearities without using Lorentz boost vector fields.
Also, it seems clear that the techniques of Sideris \cite{Si} can
handle the special case of Theorem \ref{theorem0.2} where the
semilinear terms $B^I(du)$ are not present.

We eventually choose $N = 10$ in Theorem \ref{0.2} (see \eqref{4.2} below).
The decay we obtain up until time $T_\varepsilon$ is described in
equations \eqref{Mk}-\eqref{4.4} below.

Global existence in three space dimensions has been shown for
coupled multiple speed systems satisfying various multiple-speed
versions of the so-called null condition \cite{knull}. See
Sideris-Tu \cite{Si3}, Sogge \cite{So2} for such global, multiple
speed results and further references.  These results generalize
the first global existence results of Christodoulou \cite{christ}
and  Klainerman \cite{knull}. Long time existence for multiple
speed systems in two space dimensions was studied in Kovalyov
\cite{kovalyov}.
%Global existence on the exterior of star shaped
%obstacles in three space dimensions was shown in
%\cite{KSS}, again under the assumption of a null condition.
%This extended  earlier spherically symmetric work of \cite{godin}.
%For related work outside of obstacles in higher dimensions
%see  \cite{hayashi}, \cite{tsutsumi}.

This paper is organized as follows.  In the next section we shall
prove some new pointwise $L^1\to L^\infty$ estimates for the
inhomogeneous wave equation in Minkowski space that are well
adapted to our approach of trying to mainly exploit $1/|x|$ decay
of solutions of nonlinear wave equations.  From the point of view of 
the Minkowski
space argument of Theorem \ref{theorem0.2}, this estimate is a
departure from the approach of Klainerman and Sideris \cite{KS}.
After this, we recall the weighted
space-time $L^2$ estimates from \cite{KSS2} and give
the straightforward iteration argument which
proves Theorem \ref{theorem0.2}.  We then turn to the obstacle case, obtaining
versions of the pointwise decay, weighted space-time $L^2$ estimates, and fixed time
$L^2$ estimates in the exterior of a star-shaped obstacle. As pointed out
above, the energy estimates for the boundary value problem are more involved than
their Minkowski space analogs, and in fact our estimates
involving the Euclidean rotation or scaling vector fields involve
a slight loss over their Minkowski variants.   Fortunately this
loss is not important for our goal of proving Theorem
\ref{theorem0.1}.   Finally, in   \S \ref{section:almostglobaloutside},  we combine
the decay, weighted $L^2(\R^{1+3})$, and energy estimates outside
of obstacles in an
adaptation of the proof of Theorem \ref{theorem0.2} to obtain our almost
global existence results for quasilinear wave equations outside of
star-shaped obstacles.

\newsection{Pointwise estimates in Minkowski space}

We write $\{\Omega\}=\{\Omega_{ij}\}$, where
 \begin{align}
\label{omega}
 \Omega_{ij}& =x_i\partial_j-x_j\partial_i, \quad \quad 1\le i<j\le 3,
\end{align}
 are the Euclidean $\R^3$ rotation operators. Denote by $Z$ either 
 a space-time translation or spatial rotation vector field,
 \begin{equation}
\{Z\} =\{\partial_t, \partial_j, \Omega_{ij}\}.
\label{Z}
\end{equation}
  We
 also use the scaling operator
\begin{align}
\label{scaling}
L& =t\partial_t+x\cdot\nabla_x=t\partial_t+r\partial_r.
\end{align}
Throughout the remainder of the paper we will use without explicit
mention the following fact:  if we denote by $\Gamma$ any of the vectorfields
in \eqref{Z}--\eqref{scaling}, then 
\begin{align*}
[\Gamma_i, \Gamma_j] & = \sum_{k} \mu_{ijk} \Gamma_k
\end{align*}
for certain (possibly vanishing) fixed constants $\mu_{ijk}$.

To simplify the notation, we  let
$$\square =\partial^2_t-\Delta$$
be the scalar unit-speed D'Alembertian.  We shall state most of
our estimates in terms of it, rather than the multiple speed
operator $\square_c$ in \eqref{squarem} since straightforward
scaling arguments will show that our estimates for $\square$ yield
ones for $\square_c$.

Having set up the notation, we can now state one of our main
results, which is the following variant of an
 estimate of John, Klainerman, and H\"ormander (\cite{H}, Lemma 6.6.8).

 \begin{proposition}\label{mainprop}  If $w\in C^5$ and $\square w=F$ in
 $[0,t]\times \R^3$, and the Cauchy data of $w$ vanishes
 at $t=0$, then
 \begin{equation}\label{1}
 (1+t)\,|w(t,x)| \le C\int_0^t \int_{\R^3}\sum_{|\alpha| +j \le
 3, j\le 1} |L^j Z^\alpha F(s,y)| \,\frac{dy\, ds}{|y|}
 \end{equation}
 \end{proposition}

To prove this estimate we use the following

 \begin{lemma}\label{radial}  Let $w$ be as above, and fix
 $x \in \R^3$ with $|x| = r$.   Then,
 \begin{equation}\label{rad0}|x|\, |w(t,x)|\le \frac12
 \int_0^t\int_{|r-(t-s)|}^{r+t-s}
 \sup_{|\theta|=1}|F(s,\rho\theta)|\, \rho\, d\rho\, ds.
 \end{equation}
 \end{lemma}

\begin{proof}(Lemma \ref{radial}).
This result is well known (see e.g., page 8 of Sogge \cite{S}).
Since the fundamental solution of the wave equation in $1+3$
dimensions is positive, we have that $|w|\le |W|$, where $W$ is
the solution of the inhomogeneous wave equation $\square W(t,y)=
G(t,|y|)$ and $G$ is the radial majorant of $F$,
\begin{align} \label{starone}
G(t,\rho)=\sup_{\theta \in S^2}|F(t,\rho\theta)|.
\end{align}
On the other hand, $W(t,y)$ is a spherically symmetric solution to
the wave equation in three space dimensions, hence $|y|W(t,y)$
satisfies the wave equation in one space dimension with forcing
term $|y|G(t,|y|)$,
\begin{align} \label{startwo}
|x|W(t,x)=\frac12\int_0^t\int_{|r-(t-s)|}^{r+t-s} G(s,\rho)\rho
\, d\rho \, ds. \end{align} Together, \eqref{starone}, \eqref{startwo}
yield \eqref{rad0}.
\end{proof}

\begin{proof}(Proposition \ref{mainprop}):  As in \cite{H}, we first
prove the following,
\begin{equation}\label{scaleinvariant}
 t\,|w(t,x)| \le C\int_0^t \int_{\R^3}\sum_{|\alpha|\le
 2, j\le 1} |L^j  \Omega^\alpha F(s,y)|\, \frac{dy\,ds}{|y|}.
 \end{equation}
Since the estimate \eqref{scaleinvariant} is scale invariant, it
suffices by scaling to prove the bounds for $t=1$, that is,
\begin{equation}\label{1.2}
|w(1,x)| \le C\int_0^1 \int_{\R^3}\sum_{|\alpha|\le
2, j\le 1} |L^j \Omega^\alpha F(s,y)|\, \frac{dy\,ds}{|y|}\,.
\end{equation}

Let us first prove the estimate for those $|x|>1/10$. By the
Sobolev Lemma
$$\sup_{|\theta|=1} |F(s,\rho\theta)| \le C\sum_{|\alpha|\le 2}
\int_{S^2} |(\Omega^\alpha F)(s,\rho\theta)| \, d\theta.$$
Together with \eqref{rad0} this gives,
\begin{equation}\label{11.3}
|x| \, |w(1,x)|\le C\sum_{|\alpha|\le 2}\int_0^1\int_{\R^3}
|\Omega^\alpha F(s,y)|\,\frac{dy\,ds}{|y|},\end{equation}
 which proves \eqref{1.2} when $|x|>1/10$.

It remains to consider \eqref{1.2} for a fixed  $|x|\leq 1/10$.
Since the estimate \eqref{1.2} only involves homogeneous
derivatives, and hence is preserved under cutoffs of the form
$\psi(y/|x|)$, with $\psi$ a radial bump function,
we can reduce matters to considering two cases:
 \begin{itemize}
 \item Case 1: $\supp F\subset \{(s,y): |y| \ge 2|x|\}.$
 \item Case 2: $\supp F\subset \{(s,y): |y| \le 4|x|\}.$
 \end{itemize}

 For both cases we use the formula for $w$ coming from the fundamental
 solution,
 $$w(t,x)=\frac1{4\pi}\int_{|y|<t}F(t-|y|,x-y)\, \frac{dy}{|y|}.$$

 \noindent{\bf Case 1:}  In this case $F(s,x - y)=0$ for $|y|\le
 |x|$. Hence
 $$|w(1,x)|\le \int_{|y|<1} |F(1-|y|,x-y)|\, \frac{dy}{|x-y|}.$$
 Note that $|(1-|y|,x-y)|\ge 1/4$ on the support of the integrand.
 Thus, if $\rho(s)\in C^\infty(\R)$ vanishes for $s<1/8$ and equals
 one for $s>1/4$ we have
 $$|w(1,x)|\le \int_{|y|<1} H(1-|y|,x-y)\, dy,$$
 where
 $$H(s,v)=\rho(|(s,v)|) |F(s,v)|/|v|.$$  We make the change of
 variables $\varphi(\tau,y)=\tau(1-|y|,x-y)$, where $|y|\le 1$ and
 $0<\tau<1$.
The Jacobian is
 $\tau^3(\langle x,y\rangle/|y| - 1)$.  It is bounded away from zero when
 $H(\varphi(\tau,y))\ne 0$
 since we are assuming that $|x|<1/10$, and since $H(s,v)=0$ when
 $|(s,v)|<1/8$.  Also,
\begin{multline*}\int_{|y|<1} H(1-|y|,x-y)\, dy=\int_{|y|<1 }
|H(\varphi(1,y)|\, dy
\\
\le C\iint_{|y|<1,
0<\tau<1}(|H(\varphi(\tau,y)|+|\tfrac\partial{\partial\tau}H(\varphi(\tau,y)|)
\,d\tau\, dy.
\end{multline*}
Note that $|\partial H(\varphi(\tau,y))/\partial \tau| =
|(LH)(\varphi(\tau,y)|/\tau $, and since $\tau$ is bounded from
below when $H(\varphi(\tau,y))\ne 0$, we conclude that
\begin{multline*}
|w(1,x)|\le C\iint_{0<s<1} \bigl(|H(s,y)|+|LH(s,y)|\bigr)\, dy\,ds
\\
\le C\iint_{0<s<1} \bigl(|F(s,y)|+|LF(s,y)|\bigr)\,\frac{dy\,ds}{|y|},
\end{multline*}
as desired.

\bigskip

 \noindent{\bf Case 2:}  Our assumptions here are  $F(s,y)=0$ when $|y|\ge
 4|x|$,  for some fixed $x$ with $|x|<1/10$.  In this case, we have
 $w(1,x)=w_0(1,x)$ where $w_0$ solves the inhomogeneous wave
 equation $\square w_0(t,y)=G(t,y)$, with $G(t,y)=F(t,y)$ if $t\ge
 1-5|x|$, and $G(t,y)=0$ otherwise.  By \eqref{11.3},
 \begin{align*}
 |w(1,x)|=|w_0(1,x)|&\le \frac{C}{|x|}\int^1_{1-5|x|}\int
 \sum_{|\alpha|\le 2}|\Omega^\alpha F(s,y)| \,\frac{ds\,dy}{|y|}
 \\
 &\le C\sup_{1/2<s<1}\int \sum_{|\alpha|\le 2}|\Omega^\alpha
 F(s,y)|\,\frac{dy}{|y|}.
 \end{align*}
 As in Case 1, we bound this last quantity using the fundamental
 theorem of calculus,
 \begin{align*}
 F(s,y) & = \int_0^1 \frac{d}{d\tau} G(\tau s, \tau y) \,d\tau
 \quad \frac{1}{2} \leq s \leq 1 \\
 & = \int_0^1 LG(\tau s, \tau y) \,d\tau.
 \end{align*}
 Hence we have,
\begin{align*}
|w(1,x)| & \leq C \sup_{1/2<s<1} \int_0^1 \int \sum_{|\alpha|\le
2}|\Omega^\alpha
LG(\tau s, \tau y)|\,d\tau\,\frac{dy}{|y|} \\
& \leq C \int_0^1 \int \sum_{|\alpha|\le 2}|\Omega^\alpha L F(s,y)|
\,\frac{dy}{|y|}.
\end{align*}
where, similar to Case 1 above, we've used the fact that $\left|
\frac{ \partial(\tau s, \tau y)}{\partial(\tau,y)} \right|^{-1}$
is bounded on the support of $G$.  This completes the discussion
of case 2.

Exactly as in \cite{H}, the desired bound \eqref{1} follows from
\eqref{scaleinvariant}.  More precisely, if $\text{supp}\; F
\subset \{ (s,y) \, : \, s  \geq 1 \}$, then the same is true
for $\text{supp}\; w$, and \eqref{1} follows immediately from
\eqref{scaleinvariant}.  In case $\text{supp}\; F \subset \{ (s,y)
\, : \,0 \leq s  \leq 1 \}$, then we apply the previous argument to
the function $\tilde{w}$ with $\square \tilde{w} = F(s - 2,y_1,
y_2, y_3)$. The translation introduces the usual Euclidean
derivatives and gives \eqref{1} by the preceding argument.  The
case of general forcing function $F$ follows from these
considerations and a partition of unity.
\end{proof}

\newsection{$L^2_x, L^2_{x,t}$ estimates and 
almost global existence for quasilinear equations in Minkowski space}

We now use the pointwise estimates in
Proposition \ref{mainprop} along with $L^2_tL^2_x$ estimates
exploiting $1/r$ decay of solutions of the wave equation to prove
Theorem \ref{theorem0.2}, the almost global existence theorem for
certain multiple speed systems. As we shall see,  this proof
provides a simple model for the proof of almost global existence
results in the presence of obstacles.

To do this we need to use a simple modification of an
estimate from \cite{KSS2} which
involves the scalar D'Alembertian $\square = \partial_t^2-\Delta$:

\begin{proposition}\label{prop3.1}  Suppose that $v$ solves the
wave equation $\square v=G$ on $\R_+\times \R^3$,
with Cauchy data $f\in\dot{H}^1\cap L^6(\R^3)\,,\,g\in L^2(\R^3)$ at $t=0$.
Then there is a constant $C$ so that
\begin{multline}\label{3.1}
\bigl(\ln (2+t)\bigr)^{-1/2} \|
\langle x\rangle^{-1/2}v'\|_{L^2([0,t]\times\R^3)}+
\|\langle x\rangle^{-1}v\|_{L^2_sL^6_x([0,t]\times\R^3)}
\\
\le C\,\|(f,g)\|_\holt
+C\,\int_0^t \|G(s,\cd)\|_{L^2(\R^3)}\,ds\,.
\end{multline}
\end{proposition}

Here, and in what follows, $v'$ denotes the space-time gradient of
$v$, i.e., $v'=(\partial_t v,\nabla_xv)$.

We sketch the proof of \eqref{3.1}; more details appear in
\cite{KSS2}.  The bound is achieved by considering separately two
regions of $\{(s,x): \, 0\le s\le t\}$.  Specifically, if the
norms on the left of \eqref{3.1} are taken over $\{(s,x): \, 0
\leq s \leq t\,,\, |x| \geq t\}$, the estimate follows immediately
from
$$\langle t\rangle^{-1/2}\bigl(\|v'\|_{L^2([0,t]\times
\R^3)}+\|v\|_{L^2_sL^6_x([0,t]\times \R^3)}\bigr)
\le\|(f,g)\|_\holt + \int_0^t\|G(s,\cd)\|_\ltrt \, ds
$$
which is in turn an immediate consequence of the standard fixed-time energy
estimate and Sobolev embedding. We remark that the condition $f\in L^6$
implies that $f$, hence $v(s,\cd)$ for all $s$, is the $\dot{H}^1$
limit of compactly supported functions, which allows us to bound
$\|v(s,\cd)\|_{L^6}\le C\,\|v'(s,\cd)\|_{L^2}\,.$

To establish \eqref{3.1} on the region
$\{(s,x):\, 0 \leq s \leq t\,,\, |x|< t \}$, we first show that
\begin{multline}
\label{3.2} \|v'\|_{L^2([0,t]\times\{|x|<1\})} +
\|v\|_{L^2_sL^6_x([0,t]\times\{|x|<1\})}
\\
\le C\,\|(f,g)\|_\holt+C\,\int_0^t \|G(s,\cd)\|_\ltrt\, ds\,.
\end{multline}
For the term involving $v'$ on the left, this can be shown using
the energy inequality and the sharp Huygens principle (see
\cite{KSS2} for details\footnote{In fact, the bound \eqref{3.2} is
also implicit in several previous works, going back at least to
\cite{morawetz}, \cite{morawetzstrauss}.}).
To handle the term in $v$, we note that by the Duhamel principal we may
take $G=0$. By Sobolev embedding, we can reduce
matters to showing that
$$
\|v\|_{L^2(\R\times\{|x|<1\})}
\\
\le C\,\|(f,g)\|_\holt\,.
$$
To verify this last estimate, we let $\chi(x)$ denote the
cutoff to the set $|x|<1$. Then by the Plancherel theorem, we have
$$
\|\chi v\|_{L^2_tL^2_x(\R\times\R^3)}=
\|(\hat{\chi}*\hat{v})(\tau,\xi)\|_{L^2_\tau L^2_\xi(\R\times\R^3)}
\le C\,\|\hat{v}(\tau,\xi)\|_{L^2_\tau L^1_\xi}\le 
C\,\|(f,g)\|_\holt\,,
$$
where the last inequality is seen by expressing $\hat{v}$ in terms of
$(\hat{f},\hat{g})$, and representing the $\xi$ 
integral in polar coordinates. Applying the Schwarz inequality to
the angular integral yields the desired bound.

A scaling argument applied to \eqref{3.2} yields
\begin{multline}
\label{scaled}
\|\langle x\rangle^{-1/2}v'\|_{L^2([0,t]\times\{R<|x|<2R\})}
+\|\langle x\rangle^{-1/2}v\|_{L^2_sL^6_x([0,t]\times\{R<|x|<2R\})}
\\
\le C\,\|(f,g)\|_\holt+C\,\int_0^t \|G(s,\cd)\|_\ltrt\,ds.
\end{multline}
The estimate for the first term in the left side of \eqref{3.1} on
$\{|x| \leq t  \}$ now follows by squaring the left hand side,
decomposing dyadically in $r$, using \eqref{scaled} for each
piece, and adding the resulting estimates.  One estimates the
second term in the left side of \eqref{3.1} using \eqref{scaled}
and the fact that this second term  involves the weight $\langle
x\rangle^{-1}$.  The extra weight of $\langle x
\rangle^{-\frac{1}{2}}$ allows us to sum the estimates for the
dyadic pieces with no growth in $t$.

In addition to this estimate we shall also need the standard energy estimate:

\begin{proposition}\label{prop3.2}  Let $\gamma^{IJ,ij}(t,x)$, $1\le
I,J\le N$, $0\le i,j\le 3$ be real $C^{0,1}$ functions satisfying
\begin{equation}\label{3.3}
\sum_{0\le I,J\le N}\sum_{0\le i,j\le 3}|\gamma^{IJ,ij}|< \frac{1}{2} 
\text{min}\,(c_I^2),
\quad 1\le I\le N, \, 0\le t\le T
\end{equation}
as well as
\begin{equation}\label{3.4}\int_0^T\sum_{0\le I,J\le N}\sum_{0\le i,j\le 3}
\|\nabla_{t,x}\gamma^{IJ,ij}(t,\cd)\|_{L^\infty(\R^3)} \, dt <1.
\end{equation}
Assume also that $\gamma^{IJ,ij}$ satisfies the symmetry condition
\begin{equation}\label{symm}
\gamma^{IJ,ij}=\gamma^{JI,ij}=\gamma^{IJ,ji}.\end{equation}
Then if
$$
(\partial^2_t-c_I^2\Delta) v^I=\sum_{J=1}^N\sum_{0\le i,j\le 3}
\gamma^{IJ,ij}\partial_i\partial_jv^J +F^I\,,\quad 1\le I\le N
$$
there is a constant $C$, independent of $\gamma^{IJ,ij}$, $F$, and $T$,
so that
\begin{equation}\label{3.5}
\|v'(t,\cd)\|_\ltrt \le C\,\|v'(0,\cd)\|_\ltrt +C\int_0^t\|F(s,\cd)\|_\ltrt \, ds,
\quad 0\le t\le T.
\end{equation}
\end{proposition}

We omit the standard proof of \eqref{3.5}, since analogous
estimates for  Dirichlet-wave equations will be proven in \S 5.
We observe here, though, that the energy estimate can fail in the absence of
the symmetry assumption \eqref{symm}. To see this,
consider the following nonsymmetric linear
homogeneous system on $\R \times \R^3$:
\begin{align*}
\square u & = 0 \\
\square v & = \frac{1}{4}\, \partial_t^2u\,,
\end{align*}
with $\,u\,,\,v\,,\,\partial_tv\,$ all vanishing at $t=0$,
and with $\partial_tu(0,\cd)=g$.

Then $\|(u',v')(0,\cd)\|_\ltrt = \|g\|_\ltrt,$ and the standard
energy estimate shows that
\begin{align}
\|(u',v')(t,\cd)\|^2_\ltrt &= \int_{\R^3} |g(x)|^2 dx + \frac 12
\int_0^t \int_{\R^3} v_s(s,x) u_{ss}(s,x) \, dx \, ds\,.
\label{counter}
\end{align}
Using the Fourier transform and Duhamel's principle, it is
straightforward to see that
the second term on the right hand side of \eqref{counter} is comparable to
$\|g'\|^2_\ltrt$.

We shall actually require a corollary to Proposition \ref{prop3.2}
which is based on the following commutator relations
$[(\partial^2_t-c_I^2\Delta), Z]=0$, (see \eqref{Z})
and $[(\partial^2_t-c_I^2\Delta), L]=2(\partial^2_t-c_I^2\Delta)$,
where, as above, $L$ is the scaling vector field \eqref{scaling}.

\begin{corr}\label{cor3.3}
Let $\gamma^{IJ,ij}(t,x)\in
C^\infty$ satisfy \eqref{3.3}--\eqref{symm}, and let $v$ and $F$ be
as in Proposition \ref{prop3.2}.   Then if $M=1,2,\dots$ is fixed
there is a constant $C$, independent of $\gamma^{IJ,ij}$, $F$, and
$T$, so that for $0\le t\le T$
\begin{align}\label{3.6}
\sum_{|\alpha|+m\le M}&\|L^mZ^\alpha v'(t,\cd)\|_\ltrt
\\
&\le C\sum_{|\alpha|+m\le M}\|L^mZ^\alpha v'(0,\cd)\|_\ltrt +C\int_0^t
\sum_{|\alpha|+m\le M}\|L^mZ^\alpha F(s,\cd)\|_\ltrt\, ds \notag
\\
&+ C\int_0^t\sum_{\substack{|\alpha|+m\le M \\
I,J, i,j}}\| [L^mZ^\alpha, \gamma^{IJ,ij}
\partial_i\partial_j]v^J(s,\cd)\|_\ltrt\, ds.
\notag
\end{align}
\end{corr}
We note that if we restrict $m\le 1$ on the left hand side then we
may take $m\le 1$ on the right hand side as well. We shall also
need the following consequence of the Sobolev lemma, see
Klainerman \cite{K}:

\begin{lemma}\label{lemma3.3}  Suppose that $h\in
C^\infty(\R^3)$. Then for $R>1$
$$\|h\|_{L^\infty(R/2< |x|< R)}\le
CR^{-1}\sum_{|\alpha|+|\gamma|\le 2}\|\Omega^\alpha
\partial^\gamma_x h\|_{L^2(R/4< |x|< 2R)}.
$$
\end{lemma}

To handle certain higher order commutator terms that arise in our
arguments, we will also use the following variant of an
estimate of Klainerman and Sideris (see \cite{KS}, Lemma 3.1).

\begin{lemma}\label{lemmaks}  Suppose that $1\le R\le ct/4$.  Then
for $0\le j\le 3$
\begin{align}\label{ks1}
\|\partial_jv'(t,\cd)\|_{L^2(R/2<|x|<R)}&\le
C\,(1+t)^{-1}\sum_{|\alpha|+m\le 1}\|L^mZ^\alpha
v'(t,\cd)\|_{\ltfat}
\\
&+C\,R^{-1}\bigl(\, \|v'(t,\cd)\|_{L^2(R/4<|x|<2R)}
+\|v(t,\cd)\|_{L^6(R/4<|x|<2R)}\,\bigr) \notag
\\
&+C\,\|(\partial_t^2-c^2\Delta)v\|_{L^2(R/4<|x|<2R)}\,. \notag
\end{align}
Also,
\begin{align}\label{ks2}
\|\partial_jv'(t,\cd)\|_{L^2(|x|<1)}&\le
C\,(1+t)^{-1}\sum_{|\alpha|+m\le 1}\|L^mZ^\alpha v'(t,\cd)\|_{L^2(\R^3)}
\\
&+C\,\bigl(\,
\|v'(t,\cd)\|_{L^2(|x|<2)}+\|v(t,\cd)\|_{L^6(|x|<2)}\,\bigr) \notag
\\
&+C\,\|(\partial_t^2-c^2\Delta)v\|_{L^2(|x|<2)}. \notag
\end{align}
The constant $C$ depends only on $c$.
\end{lemma}

\noindent{\bf Proof:}  By scaling we may take the wave speed $c$
to be one.  We then use the fact (see \cite{KS}, Lemma 2.3) that
for $|x|<t/2$,
$$|\partial_t v'(t,x)|+|\Delta v(t,x)|
\le C(1+t)^{-1} \!\!\!\! \sum_{|\alpha|+m\le 1}|L^mZ^\alpha
v'(t,x)| +C|(\partial^2_t-\Delta)v(t,x)|.$$ Using this we
immediately get the estimates for $j=0$.  The other cases of
\eqref{ks1} follow from the $j=0$ bound and the fact that, for
$j,k=1,2,3,$
\begin{align*}
\|&\partial_j\partial_kv(t,\cd)\|_{L^2(R/2<|x|<R)}
\\
&\le C\|\Delta v(t,\cd)\|_{L^2(R/4<|x|<2R)} +C\sum_{|\alpha|
\le 1}R^{-2+|\alpha|}\|\partial_x^\alpha v(t,\cd)\|_{L^2(R/4<|x|<2R)}
\\
&\le C\|\Delta v(t,\cd)\|_{L^2(R/4<|x|<2R)} +CR^{-1}\bigl(
\|v'(t,\cd)\|_{L^2(R/4<|x|<2R)}+\|v\|_{L^6(R/4<|x|<2R)}\bigr)\,.
\end{align*}
The inequality \eqref{ks2} follows by a similar argument.
\qed

We now use Propositions \ref{mainprop} and
\ref{prop3.1}, along with Corollary \ref{cor3.3}, to prove Theorem
\ref{theorem0.2}. We are assuming that the data $f,g \in C^{\infty}(\R^3) 
\cap L^6(\R^3)$ satisfies the
smallness condition
\begin{equation}\label{4.2}
\sum_{|\alpha|\le 10}\|\langle x\rangle^{|\alpha|}\partial_x^\alpha f'\|_\ltrt
+ \sum_{|\alpha|\le 10}
\|\langle x\rangle^{|\alpha|}\partial_x^\alpha g\|_\ltrt
\le \varepsilon,\end{equation} where $\varepsilon>0$ is small and
we aim to show that there is a solution on
$[0,T_\varepsilon]\times \R^3$, verifying
\begin{multline}\label{solution}
\sup_{0\le t\le T_\varepsilon}\Bigl(\sum_{|\alpha|+m\le 10, m\le
1}\|L^mZ^\alpha u'(t,\cd)\|_\ltrt + (1+t)\sum_{|\alpha|\le
1}\|Z^\alpha u'(t,\cd)\|_\lirt\Bigr)
\\
+ \bigl(\ln(2+T_\varepsilon)\bigr)^{-1/2}
\sum_{|\alpha|+m\le 9, m\le 1}
\|\langle x \rangle^{-\frac{1}{2}} L^mZ^\alpha
u'\|_{L^2([0,T_\varepsilon]\times \R^3)}  \le
C\,\varepsilon,\end{multline}
 where
$T_\varepsilon = \exp(\kappa/\varepsilon)$, with $\kappa>0$ being a
uniform constant.
If the initial data is $C^\infty$, and the solution
satisfies \eqref{solution},
then standard local existence theory shows that the solution is actually
$C^\infty$ on $[0,T_\varepsilon]\times \R^3$.

Set $u_{-1}=0$, and define $u_k$, $k=0,1,2,\dots$ inductively
by letting $u_k$ solve
\begin{equation}\label{4.3}
\begin{cases}
\square_{c_I} u^I_k(t,x)=B^I(u'_{k-1})+\sum_{\substack{0\le i,j,l\le 3\\
1\le J,K\le N}} B^{IJ,ij}_{Kl}
\partial_l u_{k-1}^K
\partial_{i}\partial_ju_k^J,
\\
\qquad \qquad\qquad\qquad\qquad\qquad\qquad\qquad\qquad\quad
(t,x)\in [0,T_\varepsilon]\times \R^3, \, 1\le I\le N
\\
u_k(0, \, \cdot\, )=f, \quad \partial_tu_k(0,\, \cdot \, )=g,
\end{cases}
\end{equation}
where $\square_c$ is as in \eqref{squarem}.  Let
\begin{align}\label{Mk}
M_k(T)=&\sum_{|\alpha|+m\le 9, m\le 1}
\Bigl[\bigl(\ln(2+T)\bigr)^{-1/2}\|\langle x\rangle^{-1/2}L^mZ^\alpha
u_k'\|_{L^2([0,T]\times\R^3)}
\\
&\qquad\qquad\qquad\qquad\qquad+\|\langle x\rangle^{-1}L^mZ^\alpha
u_k\|_{L^2_tL^6_x([0,T]\times\R^3)}\Bigr] \notag
\\
&+\sup_{0\le t\le T}\sum_{|\alpha|+m\le 10, m\le 1} \|L^mZ^\alpha
u_k'(t,\, \cdot\, )\|_{\ltrt} \notag
\\
&+\sup_{0\le t\le T}(1+t)\sum_{|\alpha|\le 1}\|Z^\alpha
u_k'(t,\cd)\|_\lirt \notag
\\
&=I_k(T)+II_k(T)+III_k(T)\, . \notag
\end{align}

We first observe that there is a uniform constant $C_0$ so that
$$M_0(T)\le C_0\,\varepsilon,$$
for all $T$.  This follows from
the results of section 2 and the earlier $L^2$ estimates of this section, together with
an application of the generalized Sobolev
inequalities of Klainerman to obtain the pointwise decay estimates.

We claim that if $\varepsilon<\varepsilon_0$ is sufficiently small
and if the constant $\kappa$ occurring in the definition  of
$T_\varepsilon$ is small enough, then there is a uniform constant
$C$ (which will be allowed to change from line to line throughout 
this paper) so that for all $k=1,2,3,\dots$
\begin{equation}\label{4.4}
M_k(T_\varepsilon)\le C\,\varepsilon\,.
\end{equation}
We prove this inductively.  We thus assume that the bound
holds for $k-1$ and then establish it for $k$.

We begin by applying Corollary \ref{cor3.3},  with $F=B(u_{k-1}')$ and
$$
\gamma^{IJ,ij}=\gamma^{IJ,ij}(u'_{k-1})=
\sum_{l,K}B^{IJ,ij}_{Kl}\partial_lu^K_{k-1}
\,,
$$
to estimate $II_k(T)$. Note that the hypotheses \eqref{3.3}
and \eqref{3.4} on the metric perturbation are satisfied by the
induction hypothesis if $\varepsilon$ is small and $T < T_\varepsilon$.
The symmetry hypothesis \eqref{symm} is also
valid in view of our symmetry assumption \eqref{symm0} on the
quasilinear terms. We next apply Proposition \ref{prop3.1} with
$G=\square_c L^m Z^{\alpha} u_k$ to estimate $I_k(T)$. We conclude that
\begin{multline}
\label{threeineqs}
I_k + II_k \; \le \; C_0\,\varepsilon
+\;  C\int_0^{T_\varepsilon} \sum_{\substack{|\alpha|+m\le 10\\m\le 1}}
\|L^mZ^\alpha B(u_{k-1}')(s,\cd)\|_\ltrt\,ds
\\
\;+\;
C\int_0^{T_\varepsilon}\sum_{\substack{|\alpha|+m\le 9\\m\le 1}}
\|L^mZ^\alpha \square_c u_k(s,\cd)\|_\ltrt\,ds
\\
+ \; 
C\int_0^{T_\varepsilon}\sum_{I,J,i,j}\;
\sum_{\substack{|\alpha|+m\le 10\\m\le 1}}
\|[L^mZ^\alpha,\gamma^{IJ,ij}(u'_{k-1})
\partial_i\partial_{j}]u_{k}^J(s,\cd)\|_\ltrt
 \,ds.
\end{multline}
We estimate the first integral by observing that
\begin{multline}\label{4.5}
\sum_{\substack{|\alpha|+m\le 10\\m\le 1}} |L^mZ^\alpha B(u_{k-1}')| 
\le
C\sum_{\substack{|\alpha|+m\le 9\\m\le 1}} |L^mZ^\alpha u'_{k-1}|
\sum_{\substack{|\alpha|+m\le 5\\m\le 1}} |L^mZ^\alpha u'_{k-1}|
\\
+C\,|u'_{k-1}|\sum_{\substack{|\alpha|+m\le 10\\m\le 1}} |L^mZ^\alpha u'_{k-1}|
\,.
\end{multline}
We control the contribution of the second term on the right hand side
of \eqref{4.5} to \eqref{threeineqs}
using the induction hypothesis \eqref{4.4} as follows,
\begin{multline}
\label{onestar}
\int_0^{T_\varepsilon}\|u'_{k-1}(s,\cd)||_\lirt
\sum_{\substack{|\alpha|+m\le 10\\m\le 1}}
\|L^mZ^\alpha u'_{k-1}(s,\cd)\|_\ltrt\,ds
\\
\le C\varepsilon^2\int_0^{T_\varepsilon}ds/(1+s) \le C\cdot \kappa
\cdot\varepsilon\,,
\end{multline}
where $\kappa$ is the constant appearing in \eqref{0.8}.
For the first term on the right hand side in \eqref{4.5}
we apply Lemma \ref{lemma3.3}. If we fix $s$
and $R$, we note that, for $R/2<|x|<2R,$
\begin{align}
\label{twostars}
\sum_{\substack{|\alpha|+m\le 5\\m\le 1}}\!\!
|L^mZ^\alpha u'_{k-1}(s,x)| & \le
C\,(1+R)^{-1}\!\!\!\!\sum_{\substack{|\alpha|+m\le 7\\m\le 1}}
\!\!\|L^mZ^\alpha
u'_{k-1}(s,\cd)\|_{L^2(R/4<|x|<4R)}\,.
\end{align}
We can similarly bound this factor on the set $|x| \leq 1$.  
Therefore, for each fixed $s$ we have for a given $R = 2^j$, $j \geq 0$,
\begin{multline}
\label{threestars}
\sum_{\substack{|\alpha|+m\le 9\\ m\le 1}}
\sum_{\substack{|\beta|+n\le 5\\n\le 1}}
\bigl\|\bigl(\,L^mZ^\alpha u'_{k-1}(s,\cd)\,\bigr)
\bigl(\,L^nZ^\beta u'_{k-1}(s,\cd)\,\bigr)\bigr\|_{L^2(R<|x|<2R)}  \\
\leq C 2^{-j} \sum_{\substack{|\alpha|+m\le 9\\m\le 1}}
\|L^mZ^\alpha u'_{k-1}(s,\cd)\|_{L^2(R<|x|<2R)} 
\sum_{\substack{|\alpha|+m\le 5\\m\le 1}}
\|L^mZ^\alpha u'_{k-1}(s,\cd)\|_{L^\infty(R<|x|<2R)} \\
\leq C \sum_{\substack{|\alpha|+m\le 9\\m\le 1}}
\|\langle x \rangle^{-1/2}L^mZ^\alpha u_{k-1}'(s,\cd)
\|_{L^2(R<|x|<2R)}\,,
\end{multline}
with a similar bound on the set $|x| \leq 1$, where we applied the 
Sobolev Lemma.  Summing over $R=2^j$ and using the induction hypothesis,
we conclude that
\begin{multline}
\label{fourstars}
\int_0^{T_\varepsilon} 
\sum_{\substack{|\alpha|+m\le 9\\ m\le 1}}
\sum_{\substack{|\beta|+n\le 5\\ n\le 1}}
\bigl\|\,\bigl(L^mZ^\alpha u'_{k-1}(s,\cd)\bigr)
\bigl(L^nZ^\beta u'_{k-1}(s,\cd)\bigr)\,\bigr\|_\ltrt \,ds  \\
\leq C \int_0^{T_\varepsilon} 
\sum_{\substack{|\alpha|+m\le 9\\m\le 1}}
\|\langle x \rangle^{-1/2}L^mZ^\alpha u_{k-1}'(s,\cd)\|_\ltrt^2 \\
\leq C \ln(2 + T_\varepsilon) 
\Bigl( \;\sum_{\substack{|\alpha|+m\le 9\\m\le 1}}
\bigl(\,\ln(2 + T_\varepsilon)\,\bigr)^{-\frac{1}{2}} 
\|\langle x \rangle^{-1/2}L^mZ^\alpha
u_{k-1}'\|_{L^2([0, T_\varepsilon) \times \R^3)} \Bigr)^2 \\
\leq C \cdot \kappa
\cdot\varepsilon\,.
\end{multline}
We thus have shown that
\begin{align*}
\int_0^{T_\varepsilon}\sum_{|\alpha|+m\le 10, m\le 1}
\|L^mZ^\alpha B(u_{k-1}')(s,\cd)\|_\ltrt ds & \le 
%C\varepsilon^2
%\ln(T_\varepsilon) \\
C \cdot \kappa \cdot \varepsilon\,.
\end{align*}

The second integral on the right side of \eqref{threeineqs} has a quasilinear
contribution which is bounded by
\begin{multline}
\label{plus}
\int_0^{T_\varepsilon}
\sum_{\substack{|\alpha|+m\le 9 \\m\le 1}}
\|u'_{k-1}\,(L^mZ^\alpha  u''_k(s,\cd))\|_{\ltrt}\, ds \\
+ \; \int_0^{T_\varepsilon}
\sum_{\substack{|\alpha|+m\le 5 \\ m \leq 1}}
\sum_{\substack{|\beta|+n\le 8 \\ n\le 1}}
\|(L^mZ^\alpha u'_{k-1}(s,\cd))\,
(L^nZ^\beta u''_k(s,\cd))\|_{\ltrt}\, ds \\
+ \; \int_0^{T_\varepsilon}
\sum_{\substack{|\alpha|+m\le 5 \\ m \leq 1}}
\sum_{\substack{|\beta|+n\le 8 \\ n\le 1}}
\|(L^mZ^\alpha u''_k(s,\cd))\,
(L^nZ^\beta u'_{k-1}(s,\cd))\|_{\ltrt}\, ds \\
+ \;
\int_0^{T_\varepsilon} \sum_{\substack{|\alpha|+m\le 9 \\m \leq 1}}
\|u''_k(s,\cd)\,
(L^mZ^\alpha u'_{k-1}(s,\cd))\|_{\ltrt} \, ds\,.
\end{multline}
We bound the integrand in the first integral of \eqref{plus}
by taking the first factor in $L^\infty$, the second factor in
$L^2$, and arguing as in \eqref{onestar} above to bound this
term by
\begin{align*}
 C \cdot \varepsilon \cdot
M_k(T_\epsilon)\int_0^{T_\varepsilon} \frac{1}{1+s}\, ds & \leq C
\cdot \kappa \cdot M_k(T_{\varepsilon})\,.
\end{align*}
We estimate the second and third integrals in \eqref{plus} as
before, using the generalized Sobolev bound of Lemma
\ref{lemma3.3} on the first factor. The fourth integral in
\eqref{plus} is bounded by taking the $u''_k$ factor in
$L^\infty$, and arguing as before using the induction hypothesis.
Both of these estimates yield bounds of $C \cdot \kappa \cdot
M_k(T_\varepsilon)$.

The semilinear contribution from the second integral on the right of
\eqref{threeineqs} is handled exactly as we bounded the first integral on the
right of \eqref{threeineqs}.

To estimate the third integral in \eqref{threeineqs},
which involves commutators, we begin by noting that
\begin{align*}
\sum_{\substack{|\alpha|+m\le 10 \\ m\le 1}} |\,[L^mZ^\alpha,
\gamma^{IJ,ij}(u'_{k-1})\partial_i\partial_j]\,u_k\,|
\;&\le\;
C\sum_{\substack{|\alpha|+m\le 9\\ m\le 1}}|L^m Z^\alpha
u'_{k-1}|\sum_{\substack{|\alpha|+m\le 5\\ m\le 1}}|L^mZ^\alpha
u'_k|
\\
&+\;
C\sum_{\substack{|\alpha|+m\le 10\\ m\le 1}}|L^m Z^\alpha u'_{k-1}|
\,\cdot\,|u''_k|
\\
&+\;C\sum_{\substack{|\alpha|+m\le 5\\ m\le 1}}|L^mZ^\alpha u'_{k-1}|
\sum_{\substack{|\alpha|+m\le 9\\ m\le 1}}|L^m Z^\alpha u'_{k}|
\\
&+\;C\sum_{|\alpha|\le 1}|Z^\alpha
u'_{k-1}|\sum_{\substack{|\alpha|+m\le 10\\ m\le 1}}|L^mZ^\alpha
u'_k|
\\
&+\;C\,|Lu'_{k-1}|\sum_{|\alpha|\le 9}|Z^\alpha u''_k|\,.
\end{align*}
The contribution of the first four terms to the third integral in
\eqref{threeineqs} can be controlled as in the  preceding arguments:  when
one factor appears with two or fewer $Z$ type derivatives, we take
this factor out in $L^\infty$ as in \eqref{onestar} above;
for the remaining terms we argue as in \eqref{twostars}-\eqref{fourstars}.
The last term above requires a different argument as the factor 
we would like to take in $L^\infty$ now involves the scaling vector
field $L$, which is not controlled by the term $III_{k-1}(T_\varepsilon)$ (see
\eqref{Mk}).

To estimate this last term, let $c_0=\min_I\{c_I\}$.
Then, on the region $|x|>c_0\,s/4\,,$ we can apply
Lemma \ref{lemma3.3} to obtain
$$
|Lu'_{k-1}(s,x)|\le C\,(1+s)^{-1}
\sum_{|\alpha|\le 2,m\le 1}\|L^mZ^\alpha u'_{k-1}(s,\cd)\|_\ltrt\,,
$$
and we conclude, as in  \eqref{onestar}, 
$$
\int_0^{T_\varepsilon}\sum_{|\alpha|\le 9}\|Lu'_{k-1}(s,\cd)\,Z^\alpha
u''_k(s,\cd)\|_{L^2(|x|>c_0s/4)}\,ds \le C\cdot \kappa \cdot
M_k(T_\varepsilon)\,.
$$
It remains to estimate the integrand  here on the region 
$|x|\le c_0\,s/4$. To do this, we bound the factor $Lu'_{k-1}$ in $L^\infty$ 
using Lemma \ref{lemma3.3}, then apply Lemma 
\ref{lemmaks} to $Z^\alpha u_k''$.  We obtain, for $1\le R\le c_0\,s/4\,,$
\begin{align*}
&\sum_{|\alpha|\le 9}\|Lu'_{k-1}\,Z^\alpha u''_k(s,\cd)\|_{L^2(R/2<|x|<R)}
\\
&\le C\,R^{-1}\!\!\!\!\!\sum_{|\alpha|\le 2,m\le 1} \|L^mZ^\alpha
u'_{k-1}(s,\cd)\|_\ltfat
\\
&\qquad\qquad \times
\Biggl[(1+s)^{-1}\!\!\!\!\!\!\sum_{\substack{|\alpha|+m\le 10\\m\le 1}}
\|L^mZ^\alpha u'_k(s,\cd)\|_\ltfat\\
&
\qquad\qquad\qquad\qquad\qquad\qquad
+\sum_{|\alpha|\le 9} \|Z^\alpha \square_c u_k(s,\cd)\|_\ltfat
\\
&\qquad\qquad\quad\;\;+R^{-1}\sum_{|\alpha|\le 9}
\Bigl(\;\|Z^\alpha u'_k(s,\cd)\|_{L^2(R/4<|x|<2R)}
+\|Z^\alpha u_k(s,\cd)\|_{L^6(R/4<|x|<2R)}\,\Bigr)\Biggr]\,.
\end{align*}
We can control the norms over $|x|<1$ similarly. 
After squaring this estimate and summing over dyadic values of $R$,
using extra factors of $R^{-1/2}$ to make the sums converge,
we conclude that
\begin{align}\label{endbound}
\sum_{|\alpha|\le9}&\int_0^{T_\varepsilon} \|Lu'_{k-1}Z^\alpha
u''_k(s,\cd)\|_{L^2(|x|<c_0s/4)}\, ds\\
\leq & C \int_0^{T_\varepsilon} \sum_{|\alpha|\le 2,m\le 1} 
\!\!\!\!\|L^m Z^\alpha u'_{k-1}(s,\cd) \|_{\ltrt}  
\;(1+s)^{-1} \!\!\!\!\! \sum_{\substack{|\alpha|+m\le 10\\ m\le 1}}
\|L^mZ^\alpha u'_k(s,\cd)\|_\ltrt\, ds \notag \\
&+ \int_0^{T_\varepsilon} \sum_{|\alpha| \leq 2, m\leq 1} \| L^m
Z^{\alpha} u'_{k-1}(s, \cd) \|_\ltrt 
 \times \sum_{|\alpha| \leq 9} 
\| Z^{\alpha} \square_c u_k(s, \cd) \|_{\ltrt} \, ds \notag \\
& + \int_0^{T_\varepsilon} \sum_{|\alpha|\le 2, m \leq 1} 
\| \langle
x\rangle^{-1/2}L^mZ^\alpha u'_{k-1}(s,\cd)\|_\ltrt  \notag \\
& \qquad\qquad\qquad\quad
\times \sum_{|\alpha|\le 9}
\Bigl( \|\langle x\rangle^{-1/2}Z^\alpha u'_k(s,\cd)\|_\ltrt
+\|\langle x\rangle^{-1}Z^\alpha u_k(s,\cd)\|_\lsrt\Bigr)\,ds\,. 
\notag\end{align}
The argument used in bounding the
second integral in \eqref{threeineqs} yields
$$
\int_0^{T_\varepsilon}\sum_{|\alpha|\le9} \|Z^\alpha
\square_cu_k(s,\cd)\|_\ltrt\, ds\le C \cdot \kappa
\cdot\bigl(\,\varepsilon+M_k(T_\varepsilon)\bigr)\,.
$$
Plugging this into \eqref{endbound}, and applying Cauchy-Schwarz 
to the $s$ integral of the last term on the right of \eqref{endbound},
we conclude
\begin{align*}
\sum_{|\alpha|\le9}\int_0^{T_\varepsilon}
\|L&u'_{k-1}Z^\alpha u''_k(s,\cd)\|_{L^2(|x|<c_0s/4)}\, ds
\\
& \leq \ln(2 + T_\varepsilon) \cdot M_{k-1}(T_\varepsilon)\cdot
M_k(T_\varepsilon) +
M_{k-1}(T_\varepsilon) \cdot C \cdot\kappa\cdot 
\bigl(\epsilon + M_k(T_\varepsilon)\bigr) \\
& \leq C \cdot \kappa \cdot \bigl(\varepsilon + M_k(T_\varepsilon)\bigr)\,.
\end{align*}
We have shown that
\begin{align}
\label{one} I_k + II_k & \leq C_0\,\varepsilon+ C \cdot \kappa \cdot
 \bigl(\,\varepsilon+M_k(T_\varepsilon)\bigr).
\end{align}
The final step is to show that $III_k(T_\varepsilon)$ can be controlled
in this way,
\begin{equation}\label{4.7}
\sup_{0 \leq t \leq T_\varepsilon} (1+t)  \sum_{|\alpha|\le 1}
|Z^\alpha u_k'(t,x)| \; \le \; C_0\,\varepsilon + C \cdot \kappa
\cdot \bigl(\,\varepsilon+M_k(T_\varepsilon)\bigr)\,.
\end{equation}
Together, \eqref{one} and \eqref{4.7} yield
\begin{align*}
M_k(T_\varepsilon) & \leq 3 \,C_0\, \varepsilon\,,
\end{align*}
by choosing the constant $\kappa$ sufficiently small.

It suffices then to show \eqref{4.7}.
We first note that the left hand side of \eqref{4.7} is bounded by
\begin{equation}\label{4.8}
C_0\,\varepsilon+
\int_0^t \int_{\R^3}\sum_{|\alpha| + m \le 5, m\le 1} | L^m
Z^\alpha \square_c u_k(s,y)|\;\frac{dy \, ds}{|y|}\;.
\end{equation}
This follows by Proposition \ref{mainprop}, together with the fact that
the Cauchy data of $Z^\alpha u_k'$ at $t=0$ is of size $\varepsilon$
in the appropriate norm, and hence the homogeneous
solution with the same Cauchy data satisfies the desired bounds \eqref{4.7},
by the Klainerman-Sobolev inequalities \cite{K}.

We begin by handling the integral over $|y|>1$. We note that
\begin{multline} \label{3.185}
\sum_{|\alpha| + m \le 5,m\le 1}|L^m Z^\alpha \square_c u_k(s,y)|
\\
\le C\sum_{\substack{|\alpha|+m\le 7\\ m\le 1}}|L^mZ^\alpha u'_{k-1}(s,y)|
\sum_{\substack{|\alpha|+m\le 7\\ m\le 1}}\Bigl(|L^mZ^\alpha
u'_{k-1}(s,y)|+|L^mZ^\alpha u'_{k}(s,y)|\Bigr),
\end{multline}
and conclude by the Schwarz inequality and the induction hypothesis that
\begin{align*}
 &\int_0^{T_\varepsilon}  \!\!\! \int_{|y|>1}
\sum_{|\alpha| + m \le 5, m \le 1}
| L^mZ^\alpha \square_c u_k(s,y)| \,\frac{dy\,ds}{|y|} \\
& \leq C \bigg(\sum_{\substack{|\alpha|+m\le 7\\ m\le 1}} 
\|\langle y \rangle^{-\frac{1}{2}} L^m Z^{\alpha} u_{k-1} 
\|_{L^2([0, T_\varepsilon] \times \R^3)}\bigg)^2   \\
& + C\bigg(\sum_{\substack{|\alpha|+m\le 7\\ m\le 1}} 
\| \langle y \rangle^{-\frac{1}{2}} L^m Z^{\alpha} u_{k-1} 
\|_{L^2([0, T_\varepsilon] \times \R^3)}\bigg)
\bigg( \sum_{\substack{|\alpha|+m\le 7\\ m\le 1}} 
\| \langle y \rangle^{-\frac{1}{2}} L^m Z^{\alpha} u_{k} 
\|_{L^2([0, T_\varepsilon] \times \R^3)} \bigg)\\
& \qquad\qquad\qquad\leq C \cdot \ln(2+T_\e) \cdot  
\Bigl(\, M_{k-1}^2(T_\varepsilon)+ M_{k-1}(T_\e) \cdot M_k(T_\e)\,\Bigr) \\
&\qquad\qquad\qquad 
\le C \cdot \kappa \cdot 
\bigl(\,\varepsilon + M_k(T_\varepsilon)\bigr)\\
\end{align*}
as desired.

To handle the integral over $|y|<1$, we apply the Sobolev inequality and
\eqref{3.185} to obtain
\begin{multline*}
\sum_{|\alpha| + m \le 5, m\le 1}
\sup_{|y|<1}|L^mZ^\alpha \square_c u_k(s,y)| \le C\sum_{|\alpha|+m\le 9, m\le
1}||L^mZ^\alpha u'_{k-1}(s,\cd)||_{L^2(|y|<2)}
\\
\times \sum_{|\alpha|+m\le 9,m\le 1}\Bigl(||L^mZ^\alpha
u'_{k-1}(s,\cd)||_{L^2(|y|<2)}+||L^mZ^\alpha
u'_{k}(s,\cd)||_{L^2(|y|<2)}\Bigr)\,.
\end{multline*}
Since $\frac{1}{|y|} \in L^1(\R^3)$,
\begin{align*}
&\int_0^{T_\varepsilon} \!\!\! \int_{|y| \leq 1} \sum_{|\alpha|\le
4, m \le 1} |  L^mZ^\alpha \square_c u_k(s,y)| \,\frac{dy\,ds}{|y|} \\
& \leq C\bigg( \sum_{\substack{|\alpha|+m\le 9\\ m\le 1}}
\|L^m Z^\alpha u'_{k-1}\|_{L^2([0,T_\varepsilon] \times \{|y| \leq 2 \})}
\bigg)^2 \\
& + C \bigg(
\sum_{\substack{|\alpha|+m\le 9\\m\le 1}}
\|L^m Z^\alpha u'_{k-1}\|_{L^2([0,T_\varepsilon] \times \{|y| \leq 2 \})}
\bigg)
\bigg(\sum_{\substack{|\alpha|+m\le 9\\m\le 1}}
\|L^m Z^\alpha u'_{k}\|_{L^2([0,T_\varepsilon]
\times \{|y| \leq 2 \})} \bigg) \\
&\qquad\qquad\qquad 
\le C \cdot \kappa \cdot 
\bigl(\,\varepsilon + M_k(T_\varepsilon)\bigr)\\
\end{align*}
as above. We have therefore established \eqref{4.7}.

Similar arguments show that
$$\sup_{0\le t\le
T_\varepsilon}\|u'_k(t,\cd)-u'_{k-1}(t,\cd)\|_\ltrt
\; \rightarrow 0,  \; k \to \infty.$$ We conclude that $u_k$ converges to a
solution of
\eqref{0.10} that verifies \eqref{solution} with $C=3\,C_0$.
This completes the proof of Theorem 1.2. \qed

Later we will need the following observation. If we replace the
smallness condition \eqref{4.2} by
\begin{equation}\label{small2}
\sum_{|\alpha|\le N}\|\langle x\rangle^{|\alpha|}\partial_x^\alpha f' 
\|_\ltrt +
\sum_{|\alpha|\le N-1}\|\langle x\rangle^{|\alpha|}\partial_x^\alpha g
\|_\ltrt
\le \varepsilon\,,
\end{equation}
for $N\ge 10$,  then the same argument as above gives that for
$\varepsilon>0$  small, one obtains a solution on
$[0,T_\varepsilon]\times \R^3$ verifying
\begin{multline}\label{solutioN}
\sup_{0\le t\le T_\varepsilon}\sum_{|\alpha|+m\le N, m\le 1}
\|L^mZ^\alpha u'(t,\cd)\|_\ltrt
\\
+ \bigl(\ln(2+ T_\varepsilon)\bigr)^{-1/2} \sum_{|\alpha|+m\le N-1, m\le 1}
\|L^mZ^\alpha u'\|_{L^2([0,T_\varepsilon]\times \R^3)}
\le C\,\varepsilon\,,
\end{multline} 
for this value of $N$.

\newsection{Pointwise estimates outside of star-shaped obstacles}

In this section we shall consider Dirichlet-wave equations outside
of smooth, compact, star-shaped  obstacles ${\mathcal K}\subset \R^3$.
Our main goal is to show that the solution of the inhomogeneous equation
\begin{equation}\label{5.1}
\begin{cases}
\square u(t,x)=F(t,x), \quad (t,x)\in \R_+\times \R^3\backslash {\mathcal K}
\\
u(t,x)=0, \quad x\in \partial{\mathcal K}
\\
u(t,x)=0,\quad t\le 0 
\end{cases}
\end{equation}
satisfies slightly weaker pointwise estimates than those in
Proposition \ref{mainprop}.
As before,
$\square = \partial^2_t-\Delta$ denotes the unit-speed scalar D'Alembertian,
and any of the following estimates for $\square$ extend 
to estimates for $\square_c$
after applying straightforward scaling arguments.

The pointwise estimate that we can prove is the following

\begin{theorem}\label{theorem4.1}  Suppose that ${\mathcal
K}\subset\R^3$ is a star-shaped obstacle as in \eqref{i.2}.
Then each $C^{\infty}$ solution $u$ of
\eqref{5.1} satisfies, for each $\alpha$,
\begin{multline}\label{5.0}
t\,|Z^\alpha u(t,x)| \le C\int_0^t
\int_{\R^3\backslash {\mathcal K}}
\sum_{\substack{|\beta|+j\le |\alpha|+6\\ j\le 1}} |L^j Z^\beta F(s,y)|
\frac{dy\, ds}{|y|}
\\
+C\int_0^t\sum_{\substack{|\beta|+j\le
 |\alpha|+3\\ j\le 1}} \|L^j \partial_{s,y}^\beta F(s,\cd)
\|_{L^2(\R^3\backslash {\mathcal K})}\, ds\,.
 \end{multline}
\end{theorem}

As a first step, we shall see that for any obstacle, we can reduce
things to proving decay estimates for $Z^\alpha u(t,x)$ when $x$
belongs to a fixed neighborhood of the obstacle.  Here and in what
follows, we shall assume without loss of generality that
\begin{equation}\label{5.2}{\mathcal K}\subset 
\{x\in \R^3: \, |x|<1\}.\end{equation}

\begin{lemma}\label{lemma4.2}
Suppose that $u$ as in Theorem \ref{theorem4.1} and that 
${\mathcal K}$ satisfies \eqref{5.2}.  Then
\begin{multline}\label{5.3}
(1 + t) |Z^\alpha u(t,x)|\le C\int_0^t\int_{\R^3\backslash
{\mathcal K}} \sum_{|\gamma| + j \le3, j\le 1} 
|L^j Z^{\alpha+\gamma} F(s,y)| \, \frac{dy\, ds}{|y|}
\\
+C\sup_{|y|<2, 0\le s\le t}(1+s)\bigl( \, |Z^\alpha
u'(s,y)|+|Z^\alpha u(s,y)|\bigr).
\end{multline}
\end{lemma}

\noindent{\bf Proof:}  The inequality is obvious for $|x|<2$, so we
show that there is a uniform constant $C$ such that
\begin{multline}\label{5.4}
(1+t)\sup_{|x|\ge 2}|Z^\alpha u(t,x)|\le 
C\int_0^t\int_{\R^3\backslash {\mathcal K}} \sum_{|\gamma| + j \le3, j\le 1}
|L^j Z^{\alpha+\gamma}  F(s,y)| \, \frac{dy\, ds}{|y|}
\\
+C\sup_{|y|<2, 0\le s\le t}(1+s)\bigl( \, |Z^\alpha
u'(s,y)|+|Z^\alpha u(s,y)|\bigr).
\end{multline}

For this, we fix $\rho\in C^\infty(\R)$ satisfying
$\rho(r)=1$, $r\ge 2$ and $\rho(r)=0$, $r\le1$. Then
$$w(t,x)=\rho(|x|)Z^\alpha u(t,x)$$
solves the boundaryless wave equation
$$\square w(t,x) = \rho Z^\alpha F(t,x)
-2 \rho'(|x|) \frac{x}{|x|} \cdot \nabla_x Z^\alpha u(t,x) -(\Delta \rho(|x|)
Z^\alpha u(t,x)\,,
$$
with zero initial data.  We split $w=w_0+w_1$, where $\square
w_1 =\rho Z^\alpha F$.  If we apply Proposition \ref{mainprop}, we
conclude that $(1+t)|w_1(t,x)|$ is dominated by the first term in the right
side of \eqref{5.4}, and so it suffices to show that
$(1+t)|w_0(t,x)|$
is dominated by the last term in \eqref{5.4}.  Write
$$G(t,x)=-2 \rho'(|x|) \frac{x}{|x|} \cdot \nabla_x Z^\alpha u(t,x)
-(\Delta \rho(|x|)) Z^\alpha u(t,x).$$  By Lemma \ref{radial},
\begin{align}
\label{integration}
|w_0(t,x)| &\le C \frac{1}{|x|} \int_0^t \int_{| |x| -
(t-s)|}^{  |x| + (t-s)}
\sup_{|\theta|=1}
|G(s,r\theta)|\,rdr\,ds.
\end{align}
However, $G(t,x) = 0$ for
$|x| \leq 1$ and $|x| \geq 2$.  Hence the $s$ integrand in
\eqref{integration} is nonzero only when
$$
-2 \, \leq \, |x| - (t-s) \, \leq \, 2,
$$
that is,
\begin{equation*}
(t-|x|) - 2 \; \leq \; s \; \leq \; (t-|x|) + 2.
\end{equation*}
We conclude that
\begin{align*}
|w_0(t,x)| & \leq C \frac{1}{|x|} \frac{1}{1 + |t - |x||} \sup_{
\substack{(t - |x| - 2) \leq s \leq (t - |x| + 2) \\
|y| \leq 2} } ( 1 + s ) \left( |Z^\alpha u'(s,y)|+ |Z^\alpha
u(s,y)| \right)
\end{align*}
This yields immediately  the desired bounds for $|w_0(t,x)|$ and
completes the proof of Lemma \ref{lemma4.2}. \qed

To establish decay estimates for $|x|<2$ we shall use the following
local energy estimates, which follow from the exponential
decay estimates of Lax, Morawetz, and Phillips (see \cite{LMP}, 
also \cite{MRS} for
local exponential decay outside more general obstacles).

\begin{lemma}\label{lemma4.3}  Suppose that $u\in C^\infty$
satisfies \eqref{5.1}, where ${\mathcal K}\subset \R^3$
is a star-shaped obstacle as in \eqref{5.2}. Suppose also that
$F(t,x)=0$ for $|x|>4$.  Then there is a constant $c>0$ so that
\begin{equation}\label{5.6}
\|u'(t,\cd)\|_{L^2(\R^3\backslash {\mathcal K}: \,
|x|<4)} \le C\int_0^t e^{-c(t-s)} 
\|F(s,\cd)\|_{L^2(\R^3\backslash {\mathcal K})}\, ds.
\end{equation}
Consequently, under these assumptions, if $M=0,1,2,\dots$ is fixed,
\begin{multline}\label{5.7}
\sum_{\substack{|\alpha| + j \le M\\ j\le 1}}
\|(t\partial_t)^j\partial^\alpha_{t,x} u'(t,\cd)
\|_{L^2(\R^3\backslash {\mathcal K}: \, |x|<4)} \le
C\!\!\!\!\!\sum_{\substack{|\alpha| + j \le M-1\\ j\le 1}}
\|(t\partial_t)^j\partial^\alpha_{t,x}F(t,\cd)
\|_{L^2(\R^3\backslash {\mathcal K})}
\\
+ C \int_0^te^{-\frac{c}2(t-s)}\sum_{\substack{|\alpha| + j \le M\\
j\le 1}}\|(s\partial_s)^j\partial^\alpha_{s,x} F(s,\cd)
\|_{L^2(\R^3\backslash {\mathcal K})}\,ds.
\end{multline}
\end{lemma}

\begin{proof}
The first estimate is an immediate consequence of the exponential
decay estimates of Lax and Phillips.  As for \eqref{5.7}, using induction and
elliptic regularity
(see the proof of Theorem \ref{theoreml2} below) one shows 
that for all $M = 0,1,2,\ldots,$
\begin{multline}\label{nodt}
\sum_{|\alpha|\le M}
\|\partial^\alpha_{t,x} u'(t,\cd)
\|_{L^2(\R^3\backslash {\mathcal K}: \, |x|<4)} \le
C\!\!\!\!\!\sum_{|\alpha|\le M-1}
\|\partial^\alpha_{t,x}F(t,\cd)\|_\ltrt
\\
+ C \int_0^te^{-c(t-s)}\sum_{|\alpha|\le M}
\|\partial^\alpha_{s,x} F(s,\cd)\|_\ltrt\,ds.
\end{multline}
It remains to bound
\begin{equation}
\label{lastone}
\sum_{|\alpha|\leq M - 1}\|(t\partial_t)\partial^\alpha_{t,x} u'(t,\cd)
\|_{L^2(\R^3\backslash {\mathcal K}: \, |x|<4)}.
\end{equation}
Clearly $\partial_t u$ satisfies \eqref{5.1} with
forcing term $\partial_t F$.  Apply \eqref{nodt} to 
this equation  for $\partial_t u$, 
summing on the left over
$|\alpha| \leq M-1$, and multiply
both sides of the resulting inequality by $t$ to bound 
\eqref{lastone} as in \eqref{5.7}.
\end{proof}

For later use, notice that since $L=t\partial_t+r\partial_r$,
inequality \eqref{5.7} implies that if
  $ F(s,y)=0, \, \, |y|>4$, then 
\begin{multline}\label{5.8}
\sum_{|\alpha| + j \le M, j\le 1}
\|(t\partial_t)^j\partial^\alpha_{t,x} u'(t,\cd)
\|_{L^2(\R^3\backslash {\mathcal K}: \, |x|<4)} \le
C\!\!\!\!\sum_{\substack{|\alpha|+j\le M - 1\\ j\le
1}}\|L^j\partial^\alpha_{t,x}F(t,\cd)\|_\ltrt
\\
+ C \int_0^te^{-\frac{c}2(t-s)}\sum_{\substack{|\alpha|+j\le M\\
j\le 1}}\|L^j\partial^\alpha_{s,x} F(s,\cd)\|_\ltrt \,ds. 
\end{multline}

\noindent{\bf End of proof of Theorem \ref{theorem4.1}:}  Since
the coefficients of $Z$ are bounded when $|x|<2$, it suffices to show
that if $|\beta| \leq |\alpha| + 1$ (where $\alpha$ was
fixed in the statement of the Theorem)  then
 \begin{multline}\label{5.9}
t \sup_{|x|<2}|\partial^\beta_{t,x}u(t,x)| \le
C\int_0^t\sum_{\substack{|\gamma| +j\le
  |\alpha| + 3 \\ j\le 1}} 
\|L^j \partial_{s,x}^\gamma F(s,\cd)\|_{L^2(\R^3\backslash {\mathcal K})}
\,ds\,,  \\
+C\sum_{\substack{|\gamma|+j\le |\alpha| + 4\\ j\le 1, |\mu|\le 2}}
\int_0^t\int |L^j\Omega^\mu \partial^\gamma_{s,x} F(s,y)|\,\frac{dy\,ds}{|y|}
\,.
\end{multline}

Using cutoffs for the forcing terms, we can split
things into proving \eqref{5.9} for the following two cases

\begin{itemize}
\item {\bf Case 1:}  $F(s,y)=0$ if $|y|>4$ \item {\bf Case 2:}
$F(s,y)=0$ if $|y|<3$
\end{itemize}

\noindent
For either case, we shall use the following immediate consequence of the
fundamental theorem of calculus,
$$|\,t\,\partial^\beta_{t,x}u(t,x)|\le 
\int_0^t \sum_{j\le 1}|(s\partial_s)^j \partial^\beta_{s,x}u(s,x)|\,ds\,. 
$$
We apply the Sobolev Lemma to the right side, using the fact that 
$|\beta| \leq |\alpha| + 1$,
and that Dirichlet conditions allow
us to control $u$ locally by $u'$,
to conclude that
\begin{align*} 
t\sup_{|x|<2}|\partial^\beta_{t,x}u(t,x)|
&\le C \int_0^t \sum_{|\gamma|\le |\alpha|+2,j\le 1}
\|(s\partial_s)^j \partial^\gamma_{s,x}u'(s,\cd)
\|_{L^2(\R^3\backslash {\mathcal K}: \,
|x|<4)}\,ds  \\
& \leq C \int_0^t \sum_{|\gamma| + j \le |\alpha|+3,j\le 1}
\|(s\partial_s)^j \partial^\gamma_{s,x}u'(s,\cd)
\|_{L^2(\R^3\backslash {\mathcal K}: \,|x|<4)}\,ds\,.
\end{align*}
If we are in Case 1, we can apply \eqref{5.8} to get
\eqref{5.9}.

In Case 2, we need to write $u=u_0+u_r$  where $u_0$ solves the
boundaryless wave equation $\square u_0=F$ with zero initial data.
Fix $\eta\in C^\infty_0(\R^3)$ satisfying $\eta(x)=1$,
$|x|<2$, and $\eta(x)=0$, $|x|\ge 3$.  It follows that if we set
$\tilde u=\eta u_0+u_r$ then, since $\eta F=0$,
$\tilde u$ solves the Dirichlet-wave equation
$$
\square \tilde u=G =-2\nabla_x\eta \cdot \nabla_x u_0 -(\Delta \eta)u_0
$$
with zero initial data. The forcing term $G$
vanishes unless $2\le |x|\le 4$, hence by Case 1
\begin{align*} t\sup_{|x|<2}&|\partial^\beta_{t,x}u(t,x)|
=t\sup_{|x|<2}|\partial^\beta_{t,x}\tilde u(t,x)|
\\
&\le C\int_0^t\sum_{|\gamma|+j\le
 |\alpha|+3, j\le 1} \|L^j \partial_{s,x}^\gamma G(s,\cd)
\|_{L^2(\R^3\backslash {\mathcal K})} \, ds
 \\
 &\le C\int_0^t\sum_{|\gamma|+j\le
 |\alpha|+4, j\le 1}
 \|L^j \partial_{s,x}^\gamma u_0(s,\cd)\|_{L^2(2\le
 |x|\le4)}
\, ds
 \\
 &\le C\int_0^t\sum_{|\gamma|+j\le
 |\alpha|+ 4, j\le 1} \|L^j \partial_{s,x}^\gamma u_0(s,\cd)
\|_{L^\infty(2\le
 |x|\le4)}
\, ds.
\end{align*}

To finish the argument, we apply \eqref{rad0} to 
$w=L^j\partial^\gamma_{s,x}u_0$ with $j = 0,1$. Doing so yields
\begin{align*}
\|L^j\partial^\gamma_{s,x}u_0(s,\cd)\|_{L^\infty(2\le |x|\le 4)}
&\le C\int_0^s\int_{|s-\tau-\rho|\le 4}
\sup_{|\theta|=1}|L^j\partial^\gamma_{\tau,x}F(\tau,\rho\theta)|
\, \rho \, d\rho \, d\tau
\\
&\le C\sum_{|\mu|\le 2}\int_0^s\int_{|s-\tau-\rho|\le
4}|L^j\partial^\gamma_{\tau,x}\Omega^\mu F(\tau,\rho\theta)|
\, \rho \, d\rho \, d\theta \, d\tau
\\
&= C\sum_{|\mu|\le 2}\int_0^s\int_{|s-\tau-|y|\, |\le
4}|L^j\partial^\gamma_{\tau,x}\Omega^\mu F(\tau,y)|\,\frac{dy\,d\tau}{|y|}\,.
\end{align*}

Note that the sets $\Lambda_s=\{(\tau,y):\, 0\le \tau\le s, \,
|s-\tau-|y|\, |\le 4\}$ satisfy $\Lambda_s\cap
\Lambda_{s'}=\emptyset$ if $|s-s'|>20$.  Therefore, if in the
preceding inequality we sum over $|\gamma|+j\le |\alpha| +4$, $j\le 1$ and
then integrate over $s\in [0,t]$ we conclude that \eqref{5.9} must
also hold for Case 2, which completes the proof. \qed

\newsection{Fixed time $L^2$ estimates for Euclidean vector 
fields outside obstacles}

In this section we shall work with wave equations which are small
perturbations of the standard D'Alembertian $\square$ on
$\R_+\times\R^3\backslash\mathcal{K}$. We let
$\square_\gamma$ denote the second order operator given by
\begin{equation}\label{l1}
(\square_\gamma w)^I=(\partial^2_t-c_I^2\Delta)
w^I+\sum_{J=1}^N\sum_{j,k=0}^3
\gamma^{IJ,jk}(t,x)\,\partial_j\partial_kw^J\,, \quad 1\le I\le N,
\end{equation}
where the perturbation terms $\gamma^{IJ,jk}$ satisfy the symmetry
conditions \eqref{symm}. Given $T>0$ fixed, we shall assume that
$\gamma$ is uniformly small,
\begin{equation}\label{l4}
\sum_{I,J=1}^N \sum_{j,k=0}^3
\|\gamma^{IJ,jk}(t,x)
\|_{L^{\infty}([0,T]\times{\R}^3 \backslash \mathcal{K})}
\le
\delta,
\end{equation}
and we also assume that
\begin{equation}\label{l5}
\sum_{I,J=1}^N\sum_{i,j,k=0}^3\|\partial_i
\gamma^{IJ,jk}(t,x)
\|_{L^1_tL^{\infty}_x([0,T]\times{\R}^3 \backslash \mathcal{K})}
\le C_0.
\end{equation}

Under these assumptions we shall prove $L^2$ estimates for
solutions of the inhomogeneous Dirichlet-wave equation
\begin{equation}\label{l6}
\begin{cases}
\square_\gamma w=F
\\
w|_{\partial\mathcal{K}}=0
\\
w(t,x)=0,\quad t\le 0.
\end{cases}
\end{equation}

The first estimate is the standard energy estimate:

\begin{theorem}\label{theoreml1}
Assume $w\in C^2$ satisfy \eqref{l6}, and $\gamma$ as above satisfies
the symmetry conditions \eqref{symm} as well as  
\eqref{l5} and \eqref{l4} for $\delta>0$ sufficiently small. Then
\begin{equation}\label{l7}
\|w'(t,\cd)\|_{L^2(\R^3\backslash\mathcal{K})} \le
C\int_0^t \|F(s,\cd)\|_{L^2(\R^3\backslash\mathcal{K})}\,
ds, \quad 0 \le t \le T,
\end{equation}
for a uniform constant $C$ (depending on $C_0$).
\end{theorem}

Although the result is standard, we shall present its proof
since it serves as a model for the more technical variations which
are to follow.

We first define the components of the energy-momentum vector.  For
$I=1,2,\dots,N$, we let
\begin{multline}\label{em0}
e^I_0=e^I_0(w)=(\partial_0 w^I)^2 + \sum_{k=1}^3 c^2_I (\partial_k
w^I)^2
\\
+2\sum_{J=1}^N\sum_{k=0}^3 \gamma^{IJ,0k}\partial_0w^I \partial_k
w^J - \sum_{J=1}^N \sum_{j,k=0}^3
\gamma^{IJ,jk}\partial_jw^I\partial_kw^J,
\end{multline}
and for $k=1,2,3$
\begin{equation}\label{emk}
e^I_k=e^I_k(w)=-2\,c^2_I\,\partial_0w^I\partial_kw^I+2\sum_{J=1}^N\sum_{j=0}^3
\gamma^{IJ,jk}\partial_0w^I\partial_jw^J.
\end{equation}
Then
\begin{align}\label{e0}
\partial_0e_0^I =&2\,\partial_0w^I\partial_0^2w^I + 2\sum_{k=1}^3
c^2_I\partial_kw^I\partial_0\partial_kw^I +
2\,\partial_0w^I\sum_{J=1}^N\sum_{k=0}^3\gamma^{IJ,0k}\partial_0\partial_kw^J
\\
&+2\sum_{J=1}^N\sum_{k=0}^3
\gamma^{IJ,0k}\partial^2_0w^I\partial_kw^J \notag
\\
&-\sum_{J=1}^N\sum_{j,k=0}^3 \gamma^{IJ,jk}\bigl[
\partial_0\partial_j w^I\partial_kw^J + \partial_j
w^I\partial_0\partial_kw^J \bigr]
+R_0^I, \notag
\end{align}
where
$$
R^I_0=2\sum_{J=1}^N \sum_{k=0}^3
(\partial_0\gamma^{IJ,0k})\partial_0w^I\partial_kw^J
-\sum_{J=1}^N\sum_{j,k=0}^3 (\partial_0\gamma^{IJ,jk})\partial_j
w^I\partial_k w^J.
$$
Also,
\begin{align}\label{ek}
\sum_{k=1}^3\partial_ke^I_k =&-2\,\partial_0w^Ic^2_I\Delta w^I
-2\sum_{k=1}^3c_I^2\partial_kw^I\partial_0\partial_kw^I
\\
& +2\,\partial_0w^I\sum_{J=1}^N\sum_{j=0}^3\sum_{k=1}^3
\gamma^{IJ,jk}\partial_j\partial_kw^J \notag
\\
&+2\sum_{J=1}^N\sum_{j=0}^3\sum_{k=1}^3
\gamma^{IJ,jk}\partial_0\partial_k w^I\partial_jw^J
+\sum_{k=1}^3 R^I_k, \notag
\end{align}
where
$$R^I_k=2\sum_{J=1}^N\sum_{k=1}^3
(\partial_k\gamma^{IJ,jk})\partial_0w^I\partial_jw^J.$$

Note that by the symmetry conditions \eqref{symm} if we sum the
second to last term and the third to last terms in \eqref{e0} over
$I$, we get
$$-2\sum_{I,J=1}^N \sum_{j=0}^3\sum_{k=1}^3
\gamma^{IJ,jk}\partial_0\partial_kw^I\partial_jw^J,$$ which is
$-1$ times the sum over $I$ of the second to last term of
$\eqref{ek}$.  From this, we conclude that if we set
$$e_j=e_j(w)=\sum_{I=1}^N e^I_j, \quad j=0,1,2,3,$$
and
$$R=R(w',w')=\sum_{I=1}^N \sum_{k=0}^3 R^I_k,$$
then
$$\partial_t e_0+\sum_{k=1}^3\partial_ke_k  = 
2\langle \partial_tw, \square_\gamma w\rangle
+ R(w',w'),$$ with $\langle \, \cdot \, ,\, \cdot \, \rangle$
denoting the standard inner product in $\R^N$.

If we integrate this identity over
$\R^3\backslash\mathcal{K}$ and apply the divergence
theorem, we obtain
\begin{multline}\label{l9}
\partial_t\int_{\R^3\backslash \mathcal{K}}e_0(t,x) \, dx -
\int_{\partial\mathcal{K}} \sum_{j=1}^3 e_j n_j \, d\sigma
\\
=2\int_{\R^3\backslash \mathcal{K}} \langle
\partial_tw, \square_\gamma w\rangle \, dx
+\int_{\R^3\backslash \mathcal{K}}R(w',w') \, dx\,.
\end{multline}
Here, $\vec{n}$ is the outward normal to $\mathcal{K}$, and
$d\sigma$ is surface measure on $\partial\mathcal{K}$.

Since we are assuming that $w$ solves \eqref{l6}, and hence
$\partial_t w$ vanishes on $\partial\mathcal{K}$, the integrand in the last
term in the left side of \eqref{l9} vanishes identically.  Therefore,
we have
$$
\partial_t\int_{\R^3\backslash \mathcal{K}} e_0(t,x)\, dx =
2\int_{\R^3\backslash \mathcal{K}} \langle
\partial_tw, F\rangle\,  dx
+\int_{\R^3\backslash \mathcal{K}}R(w',w')\, dx.
$$
Note that if $\delta$ in \eqref{l4} is small, then
\begin{equation}
\label{page22thing}
\bigl(2\max_I\{c^2_I, c^{-2}_I\}\bigr)^{-1}|w'(t,x)|^2\le e_0(t,x)\le
2\max_I\{c^2_I, c^{-2}_I\}|w'(t,x)|^2.
\end{equation}
This yields
\begin{multline*}
\partial_t\Bigl(\,\int_{\R^3\backslash
\mathcal{K}} e_0(t,x)\, dx\Bigr)^{1/2}
\\
\le C \|F(t,\cd)\|_{L^2(\R^3\backslash \mathcal{K})}
+C\sum_{I,J=1}^N\sum_{i,j,k=0}^3\|\partial_i
\gamma^{IJ,jk}(t,\cd)\|_{\infty}
\Bigl(\,\int_{\R^3\backslash \mathcal{K}} e_0(t,x)\,dx\Bigr)^{1/2}.
\end{multline*}
The theorem now follows from \eqref{page22thing}, \eqref{l5}, 
and Gronwall's inequality. \qed

We will also need the following estimates for $L^2$ norms of
higher order derivatives.
\begin{theorem}\label{theoreml2}  Suppose that $\gamma^{IJ,jk}\in
C^\infty([0,T]\times\R^3\backslash \mathcal{K})$ satisfy the 
symmetry conditions \eqref{symm} as well as 
\eqref{l4} and \eqref{l5} where $0<\delta<1/2$ in \eqref{l4} is
small enough so that \eqref{l7} holds.  Then if $w$ solves
\eqref{l6} and if $N=0,1,2,\dots$ there is a constant $C$,
depending on $N$, $\delta$, $\mathcal{K}$, and $C_0$, so that for $0<t<T$
\begin{align}\label{l11}
\sum_{|\alpha|\le N}\|\partial_{t,x}^\alpha
w'(t,\cd)\|_{L^2(\R^3\backslash\mathcal{K})} &\le
C\int_0^t\sum_{j\le N}\|\square_\gamma
\partial_s^jw(s,\cd)\|_{L^2(\R^3\backslash\mathcal{K})}\,ds
\\
&+C \sum_{|\alpha|\le N-1}\|\square_c\partial^\alpha_{t,x}
w(t,\cd)\|_{L^2(\R^3\backslash\mathcal{K})}\,. \notag
\end{align}
where $\square_\gamma$ and $\square_c$ are as in \eqref{l1} and
\eqref{squarem}, respectively.
\end{theorem}

\noindent{\bf Proof of Theorem \ref{theoreml2}:}  We have already
observed that \eqref{l11} holds when $N=0$, so we show that if
the estimate is valid if $N$ is replaced by $N-1$, then it must
be valid for $N$.

We first observe that, as
$\partial_tw|_{\R_+\times\partial\mathcal{K}}=0$, thus
$\sum_{|\alpha|\le
N-1}\|\partial^\alpha_{t,x}(\partial_tw)'(t,\cd)\|_\ltrt$ is dominated
by the right side of \eqref{l11}.  Hence it suffices to show
that, for $N\ge 1$,
$$\sum_{|\alpha|=N}\|\partial_{x}^\alpha
\nabla_xw(t,\cd)\|_\ltrt$$ also has this property.  But,
\begin{multline}
\label{page23stuff}
\sum_{|\alpha|=N-1}\|\Delta\partial^\alpha_{x}w(t,\cd)\|_\ltrt\\
\le C\sum_{|\alpha|=N-1}\|\partial^\alpha_{x}\partial^2_tw(t,\cd)\|_\ltrt+
C\sum_{|\alpha|=N-1}\|\square_c \partial^\alpha_{x}w(t,\cd)\|_\ltrt
\,,
\end{multline}
where $C$ depends only on the wave speeds, $c_I$. As we have observed,
the first term in the right side of \eqref{page23stuff} is
dominated by the right side of \eqref{l11}, and thus the
left side of \eqref{page23stuff} is similarly bounded. By
elliptic regularity, so is
$\sum_{|\alpha|= N}\|\partial_{x}^\alpha \nabla_xw(t,\cd)\|_\ltrt$,
which completes the proof. \qed

\newsection{Weighted $L_t^2L_x^2$ estimates for the D'Alembertian
outside star-shaped obstacles}

We shall also require $L^2_tL^2_x$ estimates for the unperturbed
inhomogeneous wave equation near the obstacle.  As in section 4, 
we consider the scalar
Dirichlet-wave equation, where 
$\square = \partial^2_t-\Delta$, 
\begin{equation}\label{l12}
\begin{cases}
\square v = G
\\
v|_{\partial\mathcal{K}}=0
\\
v(t,\cd)=0, \quad t<0.
\end{cases}
\end{equation}
Just as before, our estimates here extend to solutions of non-unit speed 
scalar wave equations after a straightforward scaling argument.  
One of the required estimates is the following,

\begin{proposition}\label{propl3}
Let $v$ be as in \eqref{l12}.  Assume also that $\mathcal{K}$ is
star-shaped and contained in $\{x\in \R^3: \, |x|<1\}$.
Then there is a uniform constant $C$ so that
\begin{equation}\label{l13}
\|v'\|_{L^2([0,t]\times\R^3\backslash\mathcal{K}: \, |x|<2)} \le
C\int_0^t \|\square v(s,\cd)\|_{\ltoo} \, ds.
\end{equation}
Additionally, if $N=1,2,3,\dots$ is fixed there is a uniform
constant $C$ so that
\begin{multline}\label{l14}
\sum_{|\alpha|\le N}\|\partial^\alpha_{s,x}v'
\|_{L^2([0,t]\times\R^3\backslash\mathcal{K}: \, |x|<2)}
\\
\le C\int_0^t \sum_{m\le N}\|\square\partial^m_sv(s,\cd)\|_\ltoo \, ds
+C\sum_{|\alpha|\le
N-1}\|\square\partial^\alpha_{s,x}v
\|_{L^2([0,t]\times\R^3\backslash \mathcal{K})}.
\end{multline}
\end{proposition}

\begin{proof} The
elliptic regularity argument used in the proof of Theorem \ref{theoreml2}
shows that
\eqref{l14} is a consequence of \eqref{l13}, so we shall
prove only \eqref{l13}.

To prove \eqref{l13}, we consider first the case when 
$G(s,y) = 0$ for $|y| > 4$.
In this case, \eqref{5.6} and the Schwarz inequality give us for $0<\tau<t$,
\begin{multline*}
\|v'(\tau,\cd)\|_{L^2(\R^3\backslash\mathcal{K}: \, |x|<2)}^2
\\
\le C\Bigl(\,\int_0^\tau e^{-c(\tau-s)}\|G(s,\cd)\|_\ltoo \, ds\Bigr)\,
\Bigl(\,\int_0^t \|G(s,\cd)\|_\ltoo \, ds\Bigr) \, .
\end{multline*}
This  implies \eqref{l14} after integrating $\tau$
from $0$ to $t$.  Note in addition that applying the Schwarz inequality to
\eqref{5.6} in a slightly different way yields,
\begin{align*}
\|v'(\tau,\cd)\|_{L^2([0,t]\times\R^3\backslash\mathcal{K}:\,|x|<2)}^2
& \le
C \int_0^t \bigg( \int_0^\tau e^{-c(\tau-s)} \|G(s,\cd)\|_\ltoo ds\bigg)^2
\, d\tau \\
& \leq C \int_0^t \int_0^\tau e^{-\frac{c}{2}(\tau-s)} \|G(s,\cd)\|_\ltoo^2
\, ds \, d\tau 
\end{align*} 
again under the
assumption that $G(s,y)=\square v(s,y)=0$, $|y|>4$.
Therefore, we also have
\begin{equation}\label{temp}
\|v'\|_{L^2(\R^3\backslash\mathcal{K}: \, |x|<2)}
\le C\,
\|G\|_{L^2( [0,t]\times\R^3\backslash\mathcal{K})},\quad
\text{if}\;\;G(s,y)=0\,,\; |y|>4.
\end{equation}

To finish, we need to show that we also have \eqref{l13} when we
assume that $G(s,y)=\square v(s,y)=0$, $|y|<3$.  For this, as in
the proof of Theorem \ref{theorem4.1} we fix $\eta\in
C^\infty(\R^3)$ satisfying $\eta(x)=1$, $|x|\le 2$, and
$\eta(x)=0$, $|x|\ge 3$.  Then if we write $v=v_0+v_r$, where $v_0$
solves the boundaryless wave equation $\square v_0=G$ with zero
initial data, it follows that $\tilde v=\eta v_0+v_r$ solves the
Dirichlet-wave equation
$$
\square \tilde v=\tilde G=
-2\nabla_x\eta\cdot\nabla_xv_0-(\Delta \eta)v_0
$$
with zero initial data, since $\eta G=0$.  Also,
$\tilde v=v$ for $|x|<2$, and $\tilde G(s,y)=0$ if $|y|>4$.  So by
\eqref{temp} we have
\begin{align*}
\|v'\|_{L^2([0,t]\times\R^3\backslash\mathcal{K}: \, |x|<2)}&=
\|\tilde v'\|_{L^2([0,t]\times\R^3\backslash\mathcal{K}: \, |x|<2)}
\\
&\le C\|\square \tilde v\|_{L^2([0,t]\times\R^3\backslash\mathcal{K})}
\\
&\le C\|v_0'\|_{L^2([0,t]\times\R^3\backslash\mathcal{K}: \, |x|<4)}
\\
 &\qquad\qquad\qquad+
C\|v_0\|_{L^2([0,t]\times\R^3\backslash\mathcal{K}: \, |x|<4)}\,.
\end{align*}
One now gets \eqref{l13} for this remaining case by applying
\eqref{3.2}, since $\square v_0=G$. 
\end{proof}

\bigskip

We also shall need $L^2_tL^2_x$ estimates involving the scaling
and Euclidean rotation vector fields.

\begin{proposition}\label{propl4}
Let $v$ and $\mathcal{K}$ be as in Proposition \ref{propl3}.  
Then if $N$ is fixed there
is a constant $C$ so that
\begin{multline}\label{l15}
\sum_{\substack{|\alpha|+m\le N \\ m\le
1}}\|L^m\partial^\alpha_{s,x}v'
\|_{L^2([0,t]\times\R^3\backslash\mathcal{K}: \, |x|<2)}
\\
\le C\int_0^t\sum_{\substack{|\alpha|+m\le N \\ m\le 1}}\|\square
L^m\partial^\alpha_{s,x}v(s,\cd)\|_\ltoo ds
+C\!\!\!\!\sum_{\substack{|\alpha|+m\le N -1\\ m\le 1}} \|\square
L^m\partial_{s,x}^\alpha
v\|_{L^2([0,t]\times\R^3\backslash \mathcal{K})}.
\end{multline}
Additionally,
\begin{align}\label{l16}
\sum_{\substack{|\alpha|+|\gamma|+m\le N \\ m\le 1}}
\|L^m\Omega^\gamma\partial^\alpha_{s,x}v'&
\|_{L^2([0,t]\times\R^3\backslash\mathcal{K}: \, |x|<2)}
\\
&\le C\int_0^t\sum_{\substack{|\alpha|+|\gamma|+m\le N \\ m\le
1}}\|\square L^m\Omega^\gamma\partial^\alpha_{s,x}v(s,\cd)\|_\ltoo \, ds
\notag
\\
&+C\!\!\!\!\sum_{\substack{|\alpha|+|\gamma|+m\le N -1\\ m\le 1}}
\|\square L^m\Omega^\gamma\partial_{s,x}^\alpha
v\|_{L^2([0,t]\times\R^3\backslash \mathcal{K})}. \notag
\end{align}
\end{proposition}

\begin{proof} We first notice that \eqref{l15}
implies \eqref{l16} since
\begin{multline*}\sum_{\substack{|\alpha|+|\gamma|+m\le N \\ m\le
1}}\|L^m\Omega^\gamma\partial^\alpha_{s,x}v'
\|_{L^2([0,t]\times\R^3\backslash\mathcal{K}: \, |x|<2)}
\\
\le C \sum_{\substack{|\alpha|+m\le N \\ m\le
1}}\|L^m\partial^\alpha_{s,x}v'
\|_{L^2([0,t]\times\R^3\backslash\mathcal{K}: \, |x|<2)}.
\end{multline*}
To prove \eqref{l15}, one repeats the proof of Proposition
\ref{propl3} using \eqref{5.8} in place of \eqref{5.6}.  
\end{proof}

As in \cite{KSS2}, one can use these estimates and the estimates
for the non-obstacle case to obtain the following.

\begin{theorem}\label{theoremw2}
Let $v$ and $\mathcal{K}$ be as in Proposition \eqref{propl3}.  
Then if $N$ is fixed there
is a constant $C$ so that
\begin{multline}\label{w6}
\bigl(\ln(2+t)\bigr)^{-1/2}\sum_{|\alpha|\le N}\|\langle
x\rangle^{-1/2}\partial^\alpha_{s,x}v'\|_{L^2(
[0,t]\times\R^3\backslash\mathcal{K})} +\sum_{|\alpha|\le
N}\|\langle x\rangle^{-1}\partial^\alpha_{s,x}v\|_{L^2_sL^6_x(
[0,t]\times\R^3\backslash\mathcal{K})}
\\
\le C\int_0^t \sum_{|\alpha|\le
N}\|\square\partial^\alpha_{s,x}v(s,\cd)\|_\ltoo \, ds
+C\sum_{|\alpha|\le N-1}\|\square\partial^\alpha_{s,x}v\|_{L^2(
[0,t]\times\R^3\backslash \mathcal{K})}.
\end{multline}
Additionally,
\begin{multline}\label{w7}
\bigl(\ln(2+t)\bigr)^{-1/2}\sum_{\substack{|\alpha|+m\le N \\ m\le
1}}\|\langle
x\rangle^{-1/2}L^m\partial^\alpha_{s,x}v'
\|_{L^2([0,t]\times\R^3\backslash\mathcal{K})}
\\\qquad\qquad\qquad\qquad\qquad
+\sum_{\substack{|\alpha|+m\le N \\ m\le 1}}\|\langle
x\rangle^{-1}L^m\partial^\alpha_{s,x}v
\|_{L^2_sL^6_x([0,t]\times\R^3\backslash\mathcal{K})}
\\
\le C\int_0^t\sum_{\substack{|\alpha|+m\le N \\ m\le 1}}\|\square
L^m\partial^\alpha_{s,x}v(s,\cd)\|_\ltoo ds
+C\!\!\!\!\!\!\!\!\sum_{\substack{|\alpha|+m\le N -1\\ m\le 1}} \|\square
L^m\partial_{s,x}^\alpha
v\|_{L^2([0,t]\times\R^3\backslash \mathcal{K})},
\end{multline}
and
\begin{align}\label{w8}
\bigl(\ln(2+t)\bigr)^{-1/2}\!\!\!\!
&\sum_{\substack{|\alpha|+m\le N \\ m\le 1}}
\|\langle
x\rangle^{-1/2}L^m Z^\alpha v'
\|_{L^2([0,t]\times\R^3\backslash\mathcal{K})}
\\
&+\sum_{\substack{|\alpha|+m\le N \\ m\le 1}}\|\langle
x\rangle^{-1}L^mZ^\alpha v\|_{L^2_sL^6_x(
[0,t]\times\R^3\backslash\mathcal{K})}
\notag
\\
&\mbox{\hspace{1in}}\le C\int_0^t\sum_{\substack{|\alpha|+m\le N \\ m\le1}}
\|\square L^m Z^\alpha v(s,\cd)\|_\ltoo\,ds
\notag
\\
&\mbox{\hspace{1.4in}}
+C\!\!\!\!\!\!\sum_{\substack{|\alpha|+m\le N -1\\ m\le 1}}
\|\square L^m Z^\alpha
v\|_{L^2([0,t]\times\R^3\backslash \mathcal{K})}. \notag
\end{align}
\end{theorem}

\begin{proof}
Let us first
handle \eqref{w6} since it is the simplest.  In view of
Proposition \ref{propl3} and Sobolev embedding it suffices to prove that
\begin{multline}\label{q1}
\bigl(\ln(2+t)\bigr)^{-1/2}\sum_{|\alpha|\le N}\|\langle
x\rangle^{-1/2}\partial^\alpha_{s,x}v'
\|_{L^2([0,t]\times\R^3\backslash\mathcal{K}: \, |x|>2)}
\\
\qquad\qquad\qquad\qquad\qquad\qquad +\sum_{|\alpha|\le N}\|\langle
x\rangle^{-1}\partial^\alpha_{s,x}v
\|_{L^2_sL^6_x([0,t]\times\R^3\backslash\mathcal{K}: \, |x|>2)}
\\
\le C\int_0^t \sum_{|\alpha|\le
N}\|\square\partial^\alpha_{s,x}v(s,\cd)\|_\ltoo \, ds
+C\sum_{|\alpha|\le N-1}\|\square\partial^\alpha_{s,x}v\|_{L^2(
[0,t]\times\R^3\backslash \mathcal{K})}.
\end{multline}

Let us estimate the first term in the left side.  For this we fix
$\beta\in C^\infty(\R^3)$ satisfying $\beta(x)=1$,
$|x|\ge 2$ and $\beta(x)=0$, $|x|\le 3/2$.  By
assumption the obstacle is contained in the set $|x|<1$. It
follows that $w=\beta v$ solves the boundaryless wave equation
$$
\square w =
\beta \square v -2\nabla_x \beta\cdot \nabla_x v - (\Delta \beta)v
$$
with zero initial data, and satisfies $w(t,x)=v(t,x)$, $|x|\ge 2$.
We split $w=w_1+w_2$, where $\square w_1=\beta \square v$, and
$\square w_2 =-2\nabla_x \beta\cdot \nabla_x v - (\Delta \beta)v$.  
Note that by \eqref{3.1} we have
$$\sum_{|\alpha|\le N}
\bigl(\ln(2+t)\bigr)^{-1/2}\|\langle x\rangle^{-1/2}\partial_{s,x}^\alpha
w_1'\|_{L^2([0,t]\times \R^3)} \le C\sum_{|\alpha|\le
N}\int_0^t \|\partial^\alpha_{s,x}\square v(s,\cd) \|_\ltoo \, ds,
$$
where we've used the fact that
$$\sum_{|\alpha|\le N}\int_0^t \|\partial^\alpha_{s,x}(\beta\square
v (s,\, \cdot\, ))\|_\ltoo \, ds \le C \sum_{|\alpha|\le N}\int_0^t
\|\partial^\alpha_{s,x}\square v(s,\cd)\|_\ltoo \, ds.$$ 
To bound the first term on the left of \eqref{q1} it 
therefore suffices to prove that
\begin{multline}\label{q2}
\bigl(\ln(2+t)\bigr)^{-1/2}\sum_{|\alpha|\le N}\|\langle
x\rangle^{-1/2}\partial_{s,x}^\alpha w_2'\|_{L^2([0,t]\times\R^3:\, |x|>2)}
\\
\le C\sum_{|\alpha|\le N}\int_0^t \|\square\partial_{s,x}^\alpha v
(s,\, \cdot\, )\|_{L^2(\R^3\backslash {\mathcal K})} \,
ds + C\sum_{|\alpha|\le N-1}\|\square\partial_{s,x}^\alpha v
\|_{L^2([0,t]\times \R^3\backslash {\mathcal K})}\,.
\end{multline}
To prove \eqref{q2} we note that $G=-2\nabla_x \beta\cdot \nabla_x
v -(\Delta \beta)v=\square w_2$ vanishes unless $1<|x|<2$.  To
use this, fix $\chi\in C^\infty_0(\R)$ satisfying
$\chi(s)=0$, $|s|>2$, and $\sum_j\chi(s-j)=1$. We then split
$G=\sum_j G_j$, where $G_j(s,x)=\chi(s-j)G(s,x)$, and let
$w_{2,j}$ be the solution of the inhomogeneous wave
equation $\square w_{2,j}=G_j$ on Minkowski space with zero initial data.
By the sharp Huygens principle, the functions $w_{2,j}$ have finite
overlap, so that we have
$|\partial_{t,x}^\alpha w_2(t,x)|^2\le C\sum_j
|\partial_{t,x}^\alpha w_{2,j}(t,x)|^2$, for some uniform constant
$C$. Therefore, by \eqref{3.1} it holds that the square of the left
side of \eqref{q2} is dominated by
\begin{align*}
\sum_{|\alpha| \leq N} &\sum_j
\Bigl(\int_0^t\|\partial_{s,x}^{\alpha}G_j(s,\, \cdot\, )\|_\ltrt ds
\Bigr)^2
\\
&\le C \sum_{|\alpha| \leq N}
\|\partial_{s,x}^{\alpha}G
\|_{L^2([0,t]\times\R^3)}^2
\\
&\le C\sum_{|\alpha|\le N}\|\partial_{s,x}^\alpha
v'\|_{L^2([0,t]\times\{ 1<|x|<2\})}^2 +
C\!\!\sum_{|\alpha|\le N}\|\partial_{s,x}^\alpha v
\|_{L^2([0,t]\times\{ 1<|x|<2\})}^2
\\
&\le C \sum_{|\alpha|\le N}\|\partial_{s,x}^\alpha
v'\|_{L^2([0,t]\times \R^3\backslash {\mathcal K}: |x|<2)}^2 .
\end{align*}
Consequently, the bound \eqref{q2} follows from \eqref{l14}.
Since the second term in \eqref{q1}
can also be handled by this argument, this completes the proof of
\eqref{w6}. Inequalities \eqref{w7} and \eqref{w8} follow by a
similar argument, using \eqref{l15} and \eqref{l16} instead of \eqref{l14}.
\end{proof}

\newsection{Fixed time $L^2$ estimates involving arbitrary differential
operators outside obstacles}

In this section we work with differential operators $P=P(t,x,D)$
which are not necessarily tangent to $\partial\mathcal{K}$, but
which satisfy other conditions to be specified. We shall
prove rather crude $L^2$ estimates for $Pw$ if $w$ solves the
inhomogeneous Dirichlet-wave equation \eqref{l6} with
$\mathcal{K}$ a star-shaped obstacle. In our applications, $P$
will be a product of powers of the Euclidean translation
vector fields and the Euclidean rotation vector fields $\{\Omega\}$,
as well as the scaling vector field, $L$.  Neither $L$ nor
the fields $\{\Omega\}$ are tangent to $\partial\mathcal{K}$ (unless
$\mathcal{K}$ is a ball, in which case the $\{\Omega\}$ vanish), and therefore,
unlike in the boundaryless case, we cannot deduce $L^2$ estimates
for $Pw$ directly from the energy estimate using
commutation properties of the vector fields and $\square_c$.  On
the other hand, the nontangential components of $L$ and $\Omega$
on $\partial\mathcal{K}$ are bounded, which leads to
estimates that we can use to prove the desired existence results.

Our basic result in this context is the following.  As remarked before,
the differential operator $P$ can be thought of as 
$P = \sum_{\substack{|\alpha|+j\le M\\ j\le 1}} L^jZ^\alpha$.

\begin{proposition}\label{propf1}
Suppose that $w$ solves \eqref{l6} where the $\gamma^{IJ,jk}$ are
as in Theorem \ref{theoreml2}.
Suppose further that there is
an integer $M$ and a constant $C_0$ so that
\begin{equation}\label{f3}
|(Pw)'(t,x)| \le C_0t \sum_{|\alpha|\le
M-1}|\partial_t\partial_{t,x}^\alpha w'(t,x)|
+C_0\sum_{|\alpha|\le M}|\partial_{t,x}^\alpha w'(t,x)|, \quad
x\in
\partial\mathcal{K}.
\end{equation}
Then, if $\square_\gamma$ and $\square_c$ are as in \eqref{l1} and
\eqref{squarem}, respectively,
\begin{align}\label{f4}
\|(Pw)'(t,\cd)\|_\ltoo &\le C\int_0^t\|\square_\gamma
Pw(s,\cd)\|_\ltoo \, ds
\\
&+C\int_0^t\sum_{\substack{|\alpha|+j\le M+1\\ j\le 1}}
\|\square_c
L^j\partial^\alpha_{s,x}w(s,\cd)\|_\ltoo \, ds \notag\\
&+C\sum_{\substack{|\alpha|+j\le M \\ j\le 1} }\|\square_c
L^j\partial^\alpha_{s,x}w\|_{L^2([0,t]\times
\R^3\backslash \mathcal{K})}\,. \notag
\end{align}
\end{proposition}

\begin{proof}
The proof is similar to that of Theorem \ref{theoreml1}, except that
here we must estimate the flux terms that arise by using the trace
inequality and the bounds \eqref{f3}.

To be more specific, we need to use the analog of \eqref{l9} where
$w$ is replaced by $Pw$.  Therefore, if we now set
$$e_j=e_j(Pw), \quad j=0,1,2,3,$$
then \eqref{l9} in our context becomes
\begin{multline}\label{f5}
\partial_t\int_{\R^3\backslash \mathcal{K}}e_0(t,x) \, dx -
\int_{\partial\mathcal{K}} \sum_{j=1}^3 e_j n_j \, d\sigma
\\
=2\int_{\R^3\backslash \mathcal{K}} \langle
\partial_tPw, \square_\gamma Pw\rangle \, dx
+\int_{\R^3\backslash \mathcal{K}}R\bigl((Pw)',(Pw)'\bigr) \,dx\,,
\end{multline}
where as before $R$ is a quadratic form whose coefficients belong
to $L^1_tL^\infty_x$.  
Therefore, if as in the proof of Theorem \ref{theoreml1}, we use
\eqref{l4} and \eqref{l5} and apply Gronwall's inequality, we
conclude that if $\delta>0$ is small enough then
\begin{align*}\|(Pw)'(t,\cd)\|_\ltoo &\le C\int_0^t\|
\square_\gamma P w(s,\cd)\|_\ltoo \, ds
\\
&+
C\Bigl(\,\int_{[0,t]\times\partial\mathcal{K}}
\bigl(\,|\partial_tPw(s,x)|^2+|\nabla_xPw(s,x)|^2\,\bigr)\,
d\sigma\Bigr)^{1/2}.
\end{align*}

Recall that we are assuming ${\mathcal K}\subset \{x\in
\R^3: \, |x|<1\}$. Therefore, by \eqref{f3} and 
a trace argument we have
\begin{equation*}\Bigl(\int_{[0,t]\times\partial\mathcal{K}}
|(Pw)'(s,x)|^2\, d\sigma\Bigr)^{1/2} \le
C\sum_{\substack{|\alpha|+j\le M+1\\ j\le 1}}
\|L^j\partial_{s,x}^\alpha w'\|_{L^2([0,t]\times
\R^3\backslash\mathcal{K}: \, |x|<2)}.
\end{equation*}
One therefore gets \eqref{f4} from \eqref{l15}.
\end{proof}

As an immediate corollary we have the following

\begin{corr}\label{corrf1}
Assume that $w$ solves \eqref{l6}.  Then if $M=1,2,\dots$
\begin{align}\label{f7}
\sum_{\substack{|\alpha|+j\le M\\ j\le 1}} \|(L^jZ^\alpha
w)'(t,\cd)\|_\ltoo
&\le C\int_0^t\sum_{\substack{|\alpha|+j\le M\\
j\le 1}} \|\square_\gamma L^jZ^\alpha w(s,\cd)\|_\ltoo \, ds
\\
&+C\int_0^t\sum_{\substack{|\alpha|+j\le M+1\\ j\le 1} }
\!\!\!\!\|\square_c L^j\partial^\alpha_{s,x}w(s,\cd)\|_\ltoo \, ds \notag
\\
&+C\sum_{\substack{|\alpha|+j\le M \\ j\le 1} }\|\square_c
L^j\partial^\alpha_{s,x}w\|_{L^2([0,t]\times
\R^3\backslash \mathcal{K})}\,. \notag
\end{align}
\end{corr}

\newsection{$L^2_x$ estimates involving only the scaling and translation
vector fields outside star-shaped obstacles}

In this section we prove $L^2_x$ estimates involving a
single occurrence of the scaling vector field $L=t\partial_t +
x\cdot \nabla_x$.  Recall that the commutator of
$L$ with $\square_c$ is $2\square_c$.
For obstacle problems, the complication arises that $L$ does not
preserve Dirichlet boundary conditions.  Because of this, unlike
in the boundaryless setting, one cannot derive $L^2$ estimates for
$Lu$ just by using energy estimates.    Fortunately, though, if
one assumes that $\partial\mathcal{K}$ is {\it star-shaped} then
in the proof of the energy estimates $L$ contributes a term with a
favorable sign, as in the classical Morawetz inequality for
star-shaped domains \cite{M}. For this reason, we can estimate
$Lu'$ in $L^2$, although there is a slight loss versus the
corresponding estimates for Minkowski space. This slight loss is
reflected especially in the third and fourth terms on the right of
\eqref{s3} below.  Unlike the corresponding terms on the right
side of \eqref{f7}, these terms involve only translation
derivatives and as such are easily handled in the nonlinear
applications to follow.

To prove the estimates of this section requires strengthening the
hypotheses on the metric perturbations $\gamma^{IJ,jk}$.  We
shall assume as before that \eqref{l5} holds, but need to
strengthen \eqref{l4} to
\begin{equation}\label{s1}
\sum_{I,J,j,k}|\gamma^{IJ,jk}(t,x)|\le \delta/(1+t),
\end{equation}
with $\delta>0$ small enough so that \eqref{l11} holds.  Under
these assumptions, we have

\begin{proposition}\label{props1}  Let $w$ solve \eqref{l6}
with $\gamma$ as in \eqref{l5}, \eqref{s1}.  
Then
\begin{align}\label{s2}
\|(Lw)'(t,\cd)\|_\ltoo &\le C\int_0^t \|\square_\gamma
Lw(s,\cd)\|_\ltoo \,ds
\\
&+ C\int_0^t\sum_{|\alpha|\le
2}\|\square_c\partial^\alpha_{s,x}w(s,\cd)\|_\ltoo \,ds
\notag
\\
&+C \sum_{|\alpha|\le
1}\|\square_c\partial^\alpha_{t,x}w\|_{L^2([0,t]\times\R^3\backslash
\mathcal{K})} . \notag
\end{align}
\end{proposition}

As a corollary of this and \eqref{l11} we have the following
useful estimate.

\begin{theorem}\label{theorems2}  Let $w$ solve \eqref{l6}
with $\gamma$ as in \eqref{l5}, \eqref{s1}.  Then if
$N=0,1,2,\dots$ is fixed
\begin{align}\label{s3}
&\sum_{\substack{|\alpha|+m\le N\\ m\le
1}}\|L^m\partial^\alpha_{t,x} w'(t,\cd)\|_\ltoo
\\
&\le C\int_0^t\sum_{\substack{|\alpha|+m\le N\\ m\le1}}\!\!\!
\|\square_\gamma L^m\partial^\alpha_{s,x}w(s,\cd)\|_\ltoo \,ds
+C\!\!\!\!\!\sum_{\substack{|\alpha|+m\le N-1\\ m\le 1}}\!\!
\|\square_c L^m
\partial^\alpha_{t,x}w(t,\cd)\|_\ltoo
 \notag
\\
&\qquad+ C\int_0^t\sum_{|\alpha|\le N+1}
\|\square_c\partial^\alpha_{s,x}w(s,\cd)\|_\ltoo\, ds +C
\sum_{|\alpha|\le N}\!
\|\square_c\partial^\alpha_{s,x}w\|_{L^2([0,t]\times\R^3\backslash
\mathcal{K})} . \notag
\end{align}
\end{theorem}

The proof that Theorem \ref{theorems2} follows from Proposition \ref{props1}
requires a simple modification of the
proof that Theorem \ref{theoreml2} follows
from Theorem \ref{theoreml1}. Precisely, we first note that, if $m=0$,
then Theorem \ref{theorems2} follows from Theorem \ref{theoreml2}.
For $m=1$, we apply induction on the number of spatial derivatives
in $\alpha.$ The elliptic regularity estimate required in this step
is that, for $N\ge 2$,
\begin{multline*}
\sum_{|\alpha|= N}\|\partial_x^\alpha Lw\|_\ltoo\le
C \sum_{|\alpha|\le N-2}
\Bigl(\,\|\partial_x^\alpha \Delta Lw\|_\ltoo
+\|\partial_x^\alpha (Lw)'\|_\ltoo\,\Bigr)
\\+
C\|(Lw)|_{\partial\mathcal{K}}\,\|_{H^{N-\frac 12}(\partial\mathcal{K})}\,.
\end{multline*}
This holds locally by standard elliptic regularity
(see e.g. \cite{gilbarg}, Theorem 8.13), and the fact that $Lw$ is
locally controlled by $(Lw)'$ and the trace of $Lw$.
Using cutoff functions one can
then reduce to the boundaryless case, where only the first term on the
right is required.

Since $(Lw)|_{\partial\mathcal{K}}=
(x\cdot\partial_x w)|_{\partial\mathcal{K}}\,,$ by the trace
theorem we have
$$
\|(Lw)|_{\partial\mathcal{K}}\,\|_{H^{N-\frac 12}(\partial\mathcal{K})}\le
C\,\sum_{|\alpha|\le N}\|\partial_x^\alpha w'\|_\ltoo\,,
$$
and the right hand side involves the estimate \eqref{theorems2}
for the case $m=0$.
It remains, then, to prove Proposition \ref{props1}.

To prove \eqref{s2}, we need to use the analog of \eqref{l9} where
$w$ there is replaced by $Lw$.  Therefore, if we set now
$$e_j=e_j(Lw), \quad j=0,1,2,3,$$
then \eqref{l9} in our context becomes
\begin{multline}\label{s4}
\partial_t\int_{\R^3\backslash \mathcal{K}}e_0(t,x) \, dx -
\int_{\partial\mathcal{K}} \sum_{j=1}^3 e_j n_j \, d\sigma
\\
=2\int_{\R^3\backslash \mathcal{K}} \langle
\partial_tLw, \square_\gamma Lw\rangle \, dx
+\int_{\R^3\backslash \mathcal{K}}R\bigl((Lw)',(Lw)'\bigr)\,dx\,,
\end{multline}
where as before $R$ is a quadratic form whose coefficients belong
to $L^1_tL^\infty_x$.

We can simplify the last term on the left hand side. We first notice
that, at points $(s,x)$ belonging to $\R_+\times\partial\mathcal{K}$,
the Dirichlet boundary conditions on $w$ give us,
$$
\partial_sLw^I=s\,\partial^2_sw^I+\partial_sw^I
+\partial_s\langle x,\nabla_x\rangle\, w^I =
\partial_s\langle x,\nabla_x\rangle\, w^I
=\langle x,\vec{n}\rangle\, \partial_{\vec{n}}\partial_sw^I\,,
$$
where
$\partial_{\vec{n}} w^I = \langle\vec{n},\nabla_x\rangle\, w^I$ denotes
differentiation with respect to the outward normal to $\mathcal{K}$.
Similarly,
$$
\sum_{j=1}^3 n_j\partial_jLw^I = s\,\partial_{\vec{n}}\partial_s w^I
+\partial_{\vec{n}}(\langle x,\nabla_x\rangle\, w^I)
$$
on
$\R_+\times\partial\mathcal{K}$.
As a consequence, we have
\begin{multline*}
-\sum_{j=1}^3e_j n_j = 2\sum_{I=1}^N
\Bigl[\langle x,\vec{n}\rangle\,c_I^2\,s\,(\partial_{\vec{n}}\partial_sw^I)^2
\\
+\langle x,\vec{n}\rangle
\,c_I^2\,\partial_{\vec{n}}
\partial_sw^I\partial_n(\langle x,\nabla_x\rangle\, w^I)
-\langle x,\vec{n}\rangle\,
\partial_{\vec{n}}\partial_sw^I
\sum_{J=1}^N\sum_{j=1}^3\sum_{k=0}^3
\gamma^{IJ,jk} n_j \partial_kLw^I\Bigr]\,.
\end{multline*}
Since we are assuming \eqref{s1}, we have
$$
-\sum_{j=1}^3e_j n_j
= 2 \sum_{I = 1}^N c_I^2 \langle x,\vec{n}\rangle\, s\,
|\partial_{\vec{n}}\partial_sw^I|^2 - Q(w'',w')\,,
$$
where
$$|Q(w'',w')|\le C
\sum_{1 \le |\alpha|\le 2} |\partial^\alpha_{s,x}w|^2 $$ for some
uniform constant $C$.
Because of this, identity \eqref{s4} yields
\begin{multline*}
\partial_t\int_{\R^3\backslash \mathcal{K}}e_0(t,x)\, dx +
\int_{\partial\mathcal{K}} 2 \sum_{I=1}^N c_I^2 \langle
x,\vec{n}\rangle\,s\,|\partial_{\vec{n}}\partial_sw^I|^2\, d\sigma
\\
= \int_{\partial\mathcal{K}} Q(w'',w')\, d\sigma
+2\int_{\R^3\backslash \mathcal{K}}
\langle\partial_tLw,\square_\gamma Lw\rangle \, dx
+\int_{\R^3\backslash \mathcal{K}}R\bigl((Lw)',(Lw)'\bigr)\, dx\,.
\end{multline*}
The second term on the left hand side is positive, since
$\langle x,\vec{n}\rangle>0$ for star-shaped $\mathcal{K}$.
Hence, we can apply Gronwall's inequality to obtain
\begin{multline*}
\|(Lw)'(t,\cd)\|_\ltoo
\le C\int_0^t\|\square_\gamma Lw(s,\cd)\|_\ltoo\, ds
\\+C\Bigl(\,\sum_{1\le|\alpha|\le2}
\int_{[0,t]\times \partial\mathcal{K}}
|\partial^\alpha_{s,x}w|^2\,d\sigma\Bigr)^{1/2}\,.
\end{multline*}
The first term on the right here is contained in the
right side of \eqref{s2}. As a result, it suffices
to show that the last term in the preceding
inequality is dominated by the other  terms in the right side of \eqref{s2}.
But
$$
\Bigl(\,\sum_{1\le|\alpha|\le2}
\int_{[0,t]\times \partial\mathcal{K}}
|\partial^\alpha_{s,x}w|^2
d\sigma\Bigr)^{1/2} \le
C\sum_{|\alpha|\le 2}
\|\partial^\alpha_{s,x} w'(s,\cd)
\|_{L^2([0,t]\times\R^3\backslash \mathcal{K}: \, |x|<2)}\,,
$$
so that \eqref{l14} yields the desired bounds for this term as
well. This completes the proof of \eqref{s2}, and hence
Proposition \ref{props1}. \qed

\newsection{Main $L^2$ estimates outside star-shaped obstacles}

We shall assume here that $\mathcal{K}\subset \R^3$ is
star-shaped. We shall also assume that the $\gamma^{IJ,jk}$
satisfy \eqref{l5} and \eqref{s1}.
Then if we combine our $L^2$ estimates we have the
following useful result.
\begin{theorem}\label{theoremm1} Let $w\in C^\infty$ solve
\eqref{l6} and vanish for $t<0$. Suppose
 also that $\mathcal{K}$ is star-shaped (see \eqref{i.2})
and the $\gamma^{IJ,jk}$ are as  in \eqref{l5}, \eqref{s1}. Then
 if $N=0,1,2,\dots$ is
fixed we have
\begin{align}\label{m1}
\sum_{|\alpha|\le N+4}&\|\partial^\alpha_{t,x}w'(t,\cd)\|_\ltoo
+\negmedspace \negmedspace \negmedspace
\sum_{\substack{|\alpha|+m\le N+2\\m\le 1}}
\negmedspace \negmedspace
\|L^m\partial^\alpha_{t,x}w'(t,\cd)\|_\ltoo 
\\
&\mbox{\hspace{2.8in}}+\negmedspace \negmedspace \negmedspace 
\sum_{\substack{|\alpha|+m\le N\\m\le 1}} 
\negmedspace \negmedspace 
\|L^mZ^\alpha w'(t,\cd)\|_\ltoo 
\notag\\
&\le C\int_0^t\Bigl(\sum_{|\alpha|\le N+4}
\|\square_\gamma\partial^\alpha_{s,x}w(s,\cd)\|_\ltoo
+\!\!\!\!\!
\sum_{\substack{|\alpha|+m\le N+2\\ m\le 1}}
\!\!\!\!\!
\|\square_\gamma L^m\partial^\alpha_{s,x}w(s,\cd)\|_\ltoo
\notag \\
&\qquad\qquad\qquad \qquad\qquad\qquad \qquad\qquad +
\sum_{\substack{|\alpha|+m\le N\\
m\le 1}}\|\square_\gamma L^mZ^\alpha w(s,\cd)\|_\ltoo \Bigr)
  ds \notag  \\
 &
+C\sum_{|\alpha|\le
N+3}\|\square_\gamma\partial^\alpha_{t,x}w(t,\cd)\|_\ltoo
+C\sum_{\substack{|\alpha|+m\le N+1\\m\le 1 }} \|\square_\gamma
L^m\partial^\alpha_{t,x}w(t,\cd)\|_\ltoo \notag \\
&+C\sum_{|\alpha|\le N+2}\|\square_c
\partial^\alpha_{s,x}w\|_{L^2([0,t]\times\R^3\backslash
\mathcal{K})} +C\sum_{\substack{|\alpha|+m\le N\\ m\le
1}}\|\square_c L^m\partial^\alpha_{s,x}w
\|_{L^2([0,t]\times\R^3\backslash\mathcal{K})} . \notag
\end{align}
\end{theorem}
 To see this, let the left side be denoted by $I+II+III$,
and let
 $RHS$
denote the right side.  We then claim that
\begin{align}
 I&\le RHS
+\sum_{I,J,j,k}\sum_{|\alpha|\le
N+3}\|\gamma^{IJ,jk}\partial_j\partial_k\partial_{t,x}^\alpha
w(t,\cd)\|_\ltoo \label{m2}
 \\
II&\le RHS +C\int_0^t\sum_{|\alpha|\le N+3}
\sum_{I,J,j,k}\|\gamma^{IJ,jk}\partial_j\partial_k
\partial^\alpha_{s,x}w(s,\cd)\|_\ltoo\, ds
 \label{m3}
 \\
&\qquad\qquad\qquad+\sum_{\substack{|\alpha|+m\le N+1\\m\le 1}}
\sum_{I,J,j,k}\|\gamma^{IJ,jk}\partial_j\partial_k
L^m\partial^\alpha_{t,x}w(t,\cd)\|_\ltoo \notag \\
III&\le RHS + C\int_0^t\sum_{\substack{|\alpha|+m\le N+1\\ m\le 1}}
\sum_{I,J,j,k}\|\gamma^{IJ,jk}\partial_j\partial_k L^m
\partial^\alpha_{s,x}w(s,\cd)\|_\ltoo\, ds.\label{m4}
\end{align}
Indeed, by
\eqref{l11}, $I$ is dominated by the first and fourth terms
in the right side
of \eqref{m1} along with the last term in \eqref{m2}.  Also,
by \eqref{s3},
$II$ is dominated by the second, fifth, first, and sixth terms in the
right side of \eqref{m1} along with the last two terms in \eqref{m3}.
Lastly, by \eqref{f7}, $III$ is dominated by the third, second,
and seventh terms in the right side of \eqref{s2},
along with the last term in \eqref{m4}.  
By inequality \eqref{l4}, if $\delta$ is sufficiently small we 
can absorb the time $t$ terms on the right hand side of \eqref{m2}
and \eqref{m3} into the left hand side of \eqref{m1}.
The inequality\eqref{m1} now follows from \eqref{l5} and
\eqref{m2}--\eqref{m4} by Gronwall's inequality. \qed

Repeating this proof and using Theorem \ref{theoremw2}
yields the following result, which will be used
in the iteration argument of the next section.

\begin{corr}\label{corollarym2} Let $w$,
$\mathcal{K}$ and $\gamma^{IJ,jk}$ be as in Theorem \ref{theoremm1}.  
Then if $N$ is fixed
\begin{align}\label{m5}
\sum_{|\alpha|\le N+4}&\|\partial^\alpha_{t,x}w'(t,\cd)\|_\ltoo +
\!\!\!\sum_{\substack{|\alpha|+m\le N+2\\ m\le
1}}\|L^m\partial^\alpha_{t,x}w'(t,\cd)\|_\ltoo
\\
&\mbox{\hspace{2.8in}}
+\!\!\!\!\!\!\!
\sum_{\substack{|\alpha|+m\le N\\ m\le 1}}\|L^mZ^\alpha w'(t,\cd)\|_\ltoo
\notag\\
&+\bigl(\ln(2+t)\bigr)^{-1/2}\Bigl(\sum_{|\alpha|\le N+3}\|\langle
x\rangle^{-1/2}\partial^\alpha_{s,x}w'
\|_{L^2([0,t]\times\R^3\backslash\mathcal{K})}
\notag \\
&\mbox{\hspace{1.4in}}
+\sum_{\substack{|\alpha|+m\le N+1\\m\le 1}}
\|\langle x\rangle^{-1/2} L^m\partial^\alpha_{s,x}w'
\|_{L^2([0,t]\times\R^3
\backslash\mathcal{K})}
\notag\\
&\mbox{\hspace{2in}}
+\sum_{\substack{|\alpha|+m\le N-1\\m\le 1}}
\|\langle x\rangle^{-1/2}L^mZ^\alpha w'
\|_{L^2([0,t]\times\R^3 \backslash\mathcal{K})}\Bigr)
\notag\\
&+\Bigl(\sum_{|\alpha|\le N+3}
\|\langle x\rangle^{-1}\partial^\alpha_{s,x}w
\|_{L^2_sL^6_x([0,t]\times\R^3\backslash\mathcal{K})}
\notag \\
&\mbox{\hspace{1.4in}} +
\sum_{\substack{|\alpha|+m\le N+1\\m\le 1}}
\|\langle x\rangle^{-1}L^m\partial^\alpha_{s,x}w
\|_{L^2_sL^6_x([0,t]\times\R^3
\backslash\mathcal{K})}\notag\\
&\mbox{\hspace{2in}}
+\sum_{\substack{|\alpha|+m\le N-1\\m\le 1}}
\|\langle x\rangle^{-1}L^mZ^\alpha w
\|_{L^2_sL^6_x([0,t]\times\R^3
\backslash\mathcal{K})}\Bigr)\notag
\\
&\!\!\!\!\!\!\!\!\!\!\!\!\!\!\!\!\!\!\!\!
\le C\int_0^t\Bigl(\sum_{|\alpha|\le N+4}
\|\square_\gamma\partial^\alpha_{s,x}w(s,\cd)\|_\ltoo+
\sum_{\substack{|\alpha|+m\le N+2\\ m\le 1}}
\|\square_\gamma L^m\partial^\alpha_{s,x}w(s,\cd)\|_\ltoo
\notag \\ 
&\mbox{\hspace{2.3in}}
+\sum_{\substack{|\alpha|+m\le N\\
m\le 1}}\|\square_\gamma L^mZ^\alpha w(s,\cd)\|_\ltoo \Bigr)
 ds \notag
\\
&+C\sum_{|\alpha|\le
N+3}\|\square_\gamma\partial^\alpha_{t,x}w(t,\cd)\|_\ltoo
+C\sum_{\substack{|\alpha|+m\le N+1\\m\le 1 }}\|\square_\gamma
L^m\partial^\alpha_{t,x}w(t,\cd)\|_\ltoo \notag
\\
&+C\sum_{|\alpha|\le N+2}\|\square_c
\partial^\alpha_{s,x}w\|_{L^2([0,t]\times\R^3
\backslash\mathcal{K})} +C\sum_{\substack{|\alpha|+m\le N\\ m\le
1}}\|\square_c L^m\partial^\alpha_{s,x}w\|_{L^2([0,t]\times
\R^3\backslash\mathcal{K})} \notag
\\
&+C\sum_{\substack{|\alpha|+m\le N-2 \\ m\le 1}}\|\square_c
L^m Z^\alpha w\|_{L^2([0,t]\times
\R^3\backslash \mathcal{K})}. \notag \end{align}

\end{corr}

\newsection{Almost global existence for quasilinear wave equations
outside of star-shaped obstacles}

\label{section:almostglobaloutside}

We conclude by showing how to adapt the proof of Theorem
\ref{theorem0.2} to establish almost global existence for
the system \eqref{0.2}. As in \cite{KSS}, \cite{KSS2}, it is convenient to
reduce the Cauchy problem \eqref{0.2} to an equivalent equation with
driving force but vanishing Cauchy data, in order to avoid dealing with
compatibility conditions for the Cauchy data.  We can then set up an
iteration argument for the new equation similar to that
used in the proof of Theorem \ref{theorem0.2}.

As in the boundaryless case, we do not need to assume that the
data has compact support.  Here, however, we have to replace the
smallness condition \eqref{4.2} by
\begin{equation}\label{10.1}
\sum_{|\alpha|\le 15}
\|\langle x\rangle^{|\alpha|}\partial_x^\alpha f\|_\ltoo
+\sum_{|\alpha|\le 14}
\|\langle x\rangle^{|\alpha|+1}\partial_x^\alpha g\|_\ltoo
\le \varepsilon\,.\end{equation}
The extra number of derivatives required is due to the loss of four
derivatives in the $L^2$ estimates for the obstacle case versus
the non-obstacle case. The extra power of $\langle x\rangle$ and 
our assumption here that we control the size of $(f,g)$ as well as their 
derivatives,  are used in the steps following \eqref{pointbound} below.

To make the reduction to an equation with zero initial data,
we first note that if the data satisfies
\eqref{10.1} with $\varepsilon>0$ small, then we can construct the
solution $u$ to the system \eqref{0.2} on
the set $(t,x) \in \{0<c\,t<|x|\}\cap \{[0,T_\varepsilon) \times \R^3\backslash\mathcal{K}\}$, where 
\begin{equation}
\label{defc}
c=5\max_I c_I\,,
\end{equation}
and that on this set the solution satisfies
\begin{equation}\label{10.2}\sup_{0\le t\le \infty}
\sum_{|\alpha|\le 15}
\|\langle x\rangle^{|\alpha|}\partial^\alpha_{t,x}u(t,\cd)
\|_{L^2(\R^3\backslash{\mathcal K}\,:\,|x|>c\,t)}
\le C_0\,\varepsilon\,.
\end{equation}
To see this, we note that by scaling the $t$ variable
we may assume that $\max_I c_I=\frac 12\,.$
The local existence results in \cite{KSS} yield a solution $u$
to \eqref{0.2} on the set 
$(t,x) \in [0,2]\times \R^3\backslash\mathcal{K}$, satisfying the bounds \eqref{10.2}.\footnote{The local
existence theorem in \cite{KSS} was stated only for diagonal
systems.  However, since the proof was based only on energy
estimates, it also applies to nonlinear systems that satisfy the
symmetry condition \eqref{symm0}, using Theorem \ref{l1}.}
To see that this solution can be extended to include all $(t,x)$ with  $0<c\,t<|x|$,
we let $R\ge 4$ and consider data $(f_R,g_R)$ supported in $R/4<|x|<4R$,
which agrees with the data $(f,g)$ on the set $R/2<|x|<2R$.
Let $u_R(t,x)$ satisfy the boundaryless equation
\begin{equation}
\label{scaledproblem}
\square_c u_R=Q(du_R,R^{-1}d^2u_R)\,,
\end{equation}
with Cauchy data $\bigl(f_R(R\cd),Rg_R(R\cd)\bigr)$.
(Recall that $Q$ is the nonlinearity appearing in 
the equation, see  \eqref{0.2}, \eqref{0.3}.)  
Because of our smallness assumption \eqref{0.7} on $(f,g)$, the solution $u_R$ of
\eqref{scaledproblem} 
exists for $0<t<1$ by standard results (see e.g. \cite{H}),
and satisfies
\begin{multline*}
\sup_{0\le t \le 1}\|u_R(t,\cd)\|_{H^{15}(\R^3)}\le
C\,
\bigl(\,\|f_R(R\,\cdot\,)\|_{H^{15}(\R^3)}
+R\,\|g_R(R\,\cdot\,)\|_{H^{14}(\R^3)}
\,\bigr)\\
\le C\,R^{-3/2}\bigg(
\sum_{|\alpha|\le 15}
\|(R\partial_x)^\alpha f_R\|_\ltrt
+R\!\!\sum_{|\alpha|\le 14}
\|(R\partial_x)^\alpha g_R\|_\ltrt\bigg)\,.
\end{multline*}
The smallness condition on $|u_R'|$ implies that the wave speeds
for the quasilinear equation \eqref{scaledproblem} are bounded above by $1$.
A domain of dependence argument shows that the solutions 
$u_R(R^{-1}t,R^{-1}x)$ restricted to $\bigl|\,|x|-R\,\bigr|<\frac{R}2-t$
agree on their overlaps, and also with the local solution $u$, yielding
the solution to \eqref{0.2}
on the desired set $(t,x) \in \{0<c\,t<|x|\}\cap \{[0,T_\varepsilon) \times 
\R^3\backslash\mathcal{K}\}.$
A partition of unity argument now yields \eqref{10.2}.

We use this partial construction of the solution $u$ to start our iteration.
Fix a cutoff function $\chi\in C^\infty(\R)$ satisfying
$\chi(s)=1$ if $s\le \frac{1}{2c}$ and $\chi(s)=0$ if $s>\frac{1}{c}$,
with $c$ as in \eqref{defc}.   Set
\begin{equation*}
u_0(t,x)=\eta(t,x)\, u(t,x)\,,\qquad
\eta(t,x)=\chi(\,|x|^{-1}t\,)\,.
\end{equation*}
Note that since $|x|$ is bounded below on the complement of $\mathcal K$,
the function $\eta(t,x)$
is smooth and homogeneous of degree 0 on $(t,x) \in [0,T_\varepsilon) 
\times (\R^3\backslash\mathcal{K})$. Also,
$$\square_c u_0=\eta Q(du,d^2u)+[\square_c,\eta]\,u\,.$$
Thus, $u$ solves $\square_c u=Q(du,d^2u)$ for $0<t<T_\varepsilon$
and $x \in {\mathbb{R}}^3 \backslash \mathcal{K}$ 
if and only if $w=u-u_0$ solves
\begin{equation}\label{10.3}
\begin{cases}
\square_c w=(1-\eta)Q\bigl(d(u_0+w), d^2(u_0+w)\bigr)
-[\square_c,\eta]u
\\
w|_{\partial \mathcal{K}}=0
\\
w(t,x)=0, \quad t\le 0
\end{cases}
\end{equation}
for $0<t<T_\varepsilon$.  We emphasize that $u_0$ has been constructed 
for all $(t,x) \in [0,T_\varepsilon) \times 
(\R^3\backslash\mathcal{K})$, and the solution $u$ has been constructed
on the support of $[\square_c,\eta]$, so that \eqref{10.3} should be viewed
as a nonlinear problem for $w$.

We shall solve \eqref{10.3} by iteration.  We set $w_0=0$,
and recursively define $w_k$ for $k=1,2,\dots$ by requiring that
\begin{equation*}
\begin{cases}
\square_c w_k = (1 - \eta)Q\bigl(d(u_0+w_{k-1}),d^2(u_0+w_{k})\bigr)
-[\square_c,\eta]u
\\
w_k|_{\partial \mathcal{K}}=0
\\
w_k(t,x)=0, \quad t\le 0\,.
\end{cases}
\end{equation*}
In place of \eqref{Mk}, we now let
\begin{align*}
M_k(T)= &\sup_{0\le t\le T} \Bigl(\sum_{|\alpha|\le 14}
\|\partial^\alpha_{t,x}w_k'(t,\cd)\|_\ltoo
+\!\!\!\!\sum_{\substack{|\alpha|+m \le 12\\ m\le 1}}
\|L^m\partial^\alpha_{t,x}w_k'(t,\cd)\|_\ltoo
\\
&\qquad+\!\!\!\!\sum_{\substack{|\alpha|+m \le 10\\ m\le 1}} 
\|L^mZ^\alpha w_k'(t,\cd)\|_\ltoo +
(1+t)\!\!\sum_{|\alpha|\le 2}\|Z^\alpha w_k(t,\cd)\|_\lioo \Bigr)
\\
&+\bigl(\ln (2+ T)\bigr)^{-1/2}
\Bigl(\sum_{\substack{|\alpha|\le 13\\}}\|\langle
x\rangle^{-1/2}
\partial^\alpha_{s,x}w_k'\|_{L^2([0,T]\times
\R^3\backslash\mathcal{K})}
\\
&\mbox{\hspace{1.4in}}+\sum_{\substack{|\alpha|+m \le 11\\ m\le 1}}
\|\langle x\rangle^{-1/2}L^m
\partial^\alpha_{s,x}w_k'\|_{L^2([0,T]\times\R^3\backslash
\mathcal{K})}
\\
&\mbox{\hspace{2in}}
+\sum_{\substack{|\alpha|+m \le 9\\
m\le 1}}\|\langle x\rangle^{-1/2}L^mZ^\alpha
w_k'\|_{L^2([0,T]\times \R^3\backslash\mathcal{K})}
\Bigr)
\\
&+\sum_{\substack{|\alpha|\le 13\\}}\|\langle
x\rangle^{-1}
\partial^\alpha_{s,x}w_k\|_{L^2_sL^6_x([0,T]\times
\R^3\backslash\mathcal{K})}
\\
&\mbox{\hspace{1.4in}} +\sum_{\substack{|\alpha|+m \le 11\\ m\le
1}}\|\langle x\rangle^{-1}L^m
\partial^\alpha_{s,x}w_k\|_{L^2_sL^6_x([0,T]\times\R^3\backslash
\mathcal{K})}
\\
&\mbox{\hspace{2in}}
+\sum_{\substack{|\alpha|+m \le 9\\
m\le 1}}\|\langle x\rangle^{-1}L^mZ^\alpha
w_k\|_{L^2_sL^6_x([0,T]\times \R^3\backslash\mathcal{K})}
\end{align*}
If we let $M_0(T)$ denote the above quantity with
$w_k$ replaced by $u_0$, then we note that \eqref{10.2} together with
Lemma \ref{lemma3.3} implies
$$
\sup_{0<T<\infty}M_0(T)\le C\,\e\,.
$$
We seek to find a constant $C_1$ so that for all $k$,
\begin{equation}\label{10.5}
M_k(T_\varepsilon)\le C_1\,\varepsilon\,,
\end{equation}
provided that
$\varepsilon<\varepsilon_0$ and provided that $\varepsilon_0$ and
the constant
$\kappa$ occurring in the definition \eqref{0.8} are
sufficiently small. To do this, we proceed inductively as in \S 3, and show
that, provided $M_{k-1}(T_\e)\le C_1\,\e$, and $\e\le\kappa$, then
\begin{equation}\label{m6}
M_k(T_\varepsilon)\le C\cdot\e+
C\cdot C_1\cdot\kappa\cdot\bigl(\,M_{k-1}(T_\e)+M_k(T_\e)\,\bigr)\,,
\end{equation}
where $C$ is a universal constant.
The bound \eqref{10.5} with $C_1=2\,C$ follows from \eqref{m6},
provided $\kappa$ is sufficiently small.

We begin by estimating the fourth term in the formula for
$M_k(T_\varepsilon)$, that is, the pointwise bounds for $Z^\alpha w_k\,.$
By Theorem \ref{theorem4.1} and the support properties of $w_k$, this is bounded by
\begin{multline}\label{pointbound}
C\int_0^{T_\e} \int_{\R^3\backslash {\mathcal K}}
\sum_{\substack{|\beta| + j \le 8\\ j\le 1}} |L^j Z^\beta \square_cw_k(s,y)|
\;\frac{dy\,ds}{|y|}
\\
+C\int_0^{T_\e}\sum_{\substack{|\beta|+j\le 5\\ j\le 1}} 
\|L^j \partial_{s,y}^\beta \square_cw_k(s,\cd)
\|_{L^2(\R^3\backslash {\mathcal K})}\,ds\,.
\end{multline}
The contribution to \eqref{pointbound} where $\square_c w_k$ is
replaced by $[\square_c,\eta]u$ is bounded by $C\,\e$.
Indeed, since this term is supported in the region 
$\frac{|y|}{2c}<s<\frac{|y|}{c}$,
\begin{multline}
\label{stuffone}
\int_0^{T_\e} \int_{\R^3\backslash {\mathcal K}}
\sum_{\substack{|\beta| + j \le 8\\ j\le 1}} |L^j Z^\beta [\square_c,\eta]\,u | 
\;\frac{dy\,ds}{|y|}  \\
\leq C \int_0^{T_\e} \langle s \rangle^{-\frac{3}{4}} \int_{\R^3\backslash {\mathcal K}}
\langle y\rangle^{-\frac{1}{4} - \frac{3}{2}}
\sum_{\substack{|\beta| + j \le 8\\ j\le 1}} \langle y \rangle^{\frac{3}{2}}
|L^j Z^{\beta} [\square_c,\eta]\,u | \, dy\, ds  \\
 \leq C \sup_{0 < t < \infty} \sum_{\substack{|\beta| + j \le 8\\ j\le 1}}
\| \langle y \rangle^2 L^j Z^{\beta}[\square_c,\eta]\,u \|_{L^2(\R^3\backslash\mathcal{K})}
\end{multline}
%\begin{equation}
%\label{stuffone}
%\sum_{|\alpha|+j\le 14}\;
%\sum_{l=0}^\infty \;
%\sup_{0<s<c\,2^{l+1}}
%\|\langle x\rangle^{\frac 32}
%L^j Z^\alpha [\square_c,\eta]\,u(s,\cd)
%\|_{L^2(\R^3\backslash\mathcal{K}\,:\,2^l<|x|<2^{l+1})}\,,
%\end{equation}
%where we use the Schwarz inequality to bound the
%first term. This in turn is bounded by
%\begin{equation}
%\label{stufftwo}
%\sup_{0<t<\infty}\sum_{|\alpha|+j\le 14}
%\|\langle x\rangle^2 L^j Z^\alpha [\square_c,\eta]\,u(t,\cd)\|_\ltoo\,,
%\end{equation}
which by \eqref{10.2}, and the homogeneity of $\eta$, is bounded by $C\,\e$.
The contribution of $[\square_c,\eta]u$ to the second term
in \eqref{pointbound} 
is bounded using  a similar argument,                  
\begin{multline}\label{stufftwo}
\int_0^{T_\e}\sum_{\substack{|\beta| + j \le 5\\ j\le 1}}
\|L^j \partial_{s,y}^\beta[\square_c,\eta]u(s, \cdot)
\|_{L^2(\R^3\backslash {\mathcal K})}\,ds\, 
\\
\leq C \int_0^{T_\e} \langle s \rangle^{-\frac{3}{2}} 
\sum_{\substack{|\beta| + j \le 5\\ j\le 1}}
\|\langle y \rangle^{\frac{3}{2}}
L^j \partial_{s,y}^\beta[\square_c,\eta]u(s, \cdot)
\|_{L^2(\R^3\backslash {\mathcal K})}\,ds\,  
\end{multline}
which is bounded by $C\, \e$ as argued for \eqref{stuffone}.

The  contribution to the expression \eqref{pointbound} in which
we replace the term $\square_c w_k$ by the term
\mbox{$(1-\eta)Q\bigl(d(u_0+w_{k-1}),d^2(u_0+w_k)\bigr)$}
can be bounded 
by
\begin{multline}
\label{thebound}
C\cdot\Bigl(\,\varepsilon+
\sum_{\substack{|\alpha|+m\le 9\\m\le 1}}
\|\langle x\rangle^{-\frac 12}L^mZ^\alpha w'_{k-1}
\|_{L^2([0,T_\e]\times\R^3\backslash{\mathcal K})}\Bigr)\times\\
\Bigl(\,\varepsilon+\!\!\!\!
\sum_{\substack{|\alpha|+m\le 9\\m\le 1}}\!\!\!\!
\|\langle x\rangle^{-\frac 12}L^mZ^\alpha w'_{k-1}
\|_{L^2([0,T_\e]\times\R^3\backslash{\mathcal K})}+
\!\!\!\!\sum_{\substack{|\alpha|+m\le 9\\m\le 1}}\!\!\!\!
\|\langle x\rangle^{-\frac 12}L^mZ^\alpha w'_k
\|_{L^2([0,T_\e]\times\R^3\backslash{\mathcal K})}\Bigr)\,.
\end{multline}
For example, a typical term
involving $u_0$ coming from the first term in \eqref{pointbound} 
is handled as follows
\begin{multline}
\label{typicaluzeroterm}
\int_0^{T_\e} \int_{\R^3\backslash {\mathcal K}}
\sum_{\substack{|\beta|+j\le 8\\j\le 1}} 
(1-\eta)\,|L^j Z^{\beta} (u_0'\,w_k '')| \frac{dy\,ds}{|y|} \\
\leq C \int_0^{T_\e} 
\langle s\rangle^{-\frac 32}\!\!\!\!
\sum_{\substack{|\beta| + j \le 9\\ j\le 1}}\!\!
\|\langle y \rangle L^j  Z^{\beta} u_0'(s, \cd) 
\|_{L^2(\R^3\backslash{\mathcal K})}
\sum_{\substack{|\beta| + j \le 9\\ j\le 1}}\!\!\!
\|\langle y \rangle^{-\frac 12} L^j  Z^{\beta} w_k'(s, \cd) 
\|_{L^2(\R^3\backslash{\mathcal K})} \,ds \\
\leq C
\sum_{\substack{|\beta| + j \le 9\\ j\le 1}}
\|\langle s \rangle^{-\frac 32}\langle y \rangle L^j  Z^{\beta} u_0' 
\|^2_{L^2([0,T_\varepsilon)\times\R^3\backslash{\mathcal K})}
\sum_{\substack{|\beta| + j \le 9\\ j\le 1}}
\|\langle y \rangle^{-\frac 12} L^j  Z^{\beta} w_k' 
\|_{L^2([0, T_\varepsilon) \times \R^3\backslash{\mathcal K})}
\\
\leq C\,\varepsilon \sum_{\substack{|\beta| + j \le 9\\ j\le 1}} 
\| \langle y \rangle^{-\frac 12} L^j  Z^{\beta} w_k' 
\|_{L^2([0, T_\varepsilon) \times \R^3\backslash{\mathcal K})}
\end{multline}
again using \eqref{10.2}.
Arguing as in \eqref{typicaluzeroterm}, \eqref{stuffone} 
and \eqref{stufftwo}, one
easily checks the bound \eqref{thebound} for the other 
terms in \eqref{pointbound} involving $u_0$.  To bound
the contributions to \eqref{pointbound} involving only $w_{k-1}, w_k$, we
apply the Schwarz inequality to handle the first term in \eqref{pointbound},
 and for the second term in \eqref{pointbound}, the bound \eqref{thebound}
follows by applying
Lemma \ref{lemma3.3} in the manner used to bound the first term
on the right side of  \eqref{4.5} (see \eqref{twostars} --\eqref{fourstars}).
We thus have the bound
\begin{multline*}
C\,\Bigl(\,\varepsilon+\ln(T_\e)^{\frac 12}\,M_{k-1}(T_\e)\Bigr)\times
\Bigl(\,\varepsilon+\ln(T_\e)^{\frac 12}\,M_{k-1}(T_\e)
+\ln(T_\e)^{\frac 12}\,M_k(T_\e)
\Bigr)\\
\le
C\cdot\e+C\cdot C_1\cdot\kappa\cdot\bigl(\,M_{k-1}(T_\e)+M_k(T_\e)\,\bigr)\,.
\end{multline*}
Thus the fourth term in the definition of $M_k(T_\e)$ satisfies the bounds
\eqref{m6}.

All other terms in $M_k(t)$ occur in the left hand side of \eqref{m5},
taking $N=10$, $w = w_k$, $t=T_\varepsilon$,
and letting $\gamma$ be defined by
\begin{align}\label{writeitout}
(\square_\gamma w_k)^I &\equiv
\square_{c_I}w_k^I-
(1 - \eta)\sum_{\substack{0\le i,j,l\le 3\\ 1\le J,K\le N} }B^{IJ,ij}_{K,l}
(\partial_l u_0^K+\partial_l w_{k-1}^K)\,\partial_i\partial_j  w_k^J\\
&=(1 - \eta) B^I\bigl(d(u_0+w_{k-1})\bigr)-[\square_{c_I},\eta]u^I
\notag\\
&\mbox{\hspace{1in}}
+(1 - \eta)\sum_{\substack{0\le i,j,l\le 3\\ 1\le J,K\le N} }B^{IJ,ij}_{K,l}
(\partial_l u_0^K+\partial_l w_{k-1}^K)\,\partial_i\partial_j u_0^J 
\notag\,. 
\end{align}
Hence, we need to show that
each term on the right of \eqref{m5}, with these values for 
$w_k, \gamma,t$ and $N$,
can be dominated by the right hand side of \eqref{m6}.

We first estimate
\begin{multline}\label{firstthreeterms}
\int_0^t\Bigl(\sum_{|\alpha|\le 14}
\|\square_\gamma\partial^\alpha_{s,x}w_k(s,\cd)\|_\ltoo+
\sum_{\substack{|\alpha|+m\le 12\\ m\le 1}}
\|\square_\gamma L^m\partial^\alpha_{s,x}w_k(s,\cd)\|_\ltoo
\\ 
\mbox{\hspace{2.3in}}
+\sum_{\substack{|\alpha|+m\le 10\\m\le 1}}
\|\square_\gamma L^mZ^\alpha w_k(s,\cd)\|_\ltoo \Bigr)\,ds \,.
\end{multline}
Consider the third term in \eqref{firstthreeterms}, which is
clearly bounded by 
\begin{equation}
\label{writeouttrivial}
\sum_{\substack{|\alpha|+m\le 10\\m\le 1}}
\int_0^t\Bigl(\| L^mZ^\alpha \square_\gamma  w_k(s,\cd)\|_\ltoo
+ \| \, [\square_\gamma, L^m Z^\alpha] w_k(s, \cdot) \|_\ltoo \Bigr)\,ds \,.
\end{equation}
The contribution to  $\square_\gamma w_k$ in the first term of \eqref{writeouttrivial}
coming from $[\square_c,\eta]u$ is
handled as in \eqref{stufftwo} above.  
To handle the contribution here from 
$(1 - \eta)B^I(d(u_0 + w_{k-1}))$, 
we bound those terms involving $u_0$ using \eqref{10.2} 
arguing similar to \eqref{stuffone}, \eqref{stufftwo},
and \eqref{typicaluzeroterm} above.
In the same way, one bounds the contribution 
of the last term on the right side of 
\eqref{writeitout} to the first term in \eqref{writeouttrivial}.
   The contributions
which remain to be bounded from both terms of \eqref{writeouttrivial}
are identical
to the first and third integrals on the right hand side of
\eqref{threeineqs}, with $u'_{k-1}, u'_k$ there replaced by 
$w'_{k-1}, w'_k$, respectively. These terms can be estimated
in an identical manner
to that section, with the following remark in mind.
Lemma \ref{lemma3.3} holds on our exterior domain 
for a general function $h$ (and the same proof applies without modification),
but the analogue of estimate \eqref{ks2} of Lemma \ref{lemmaks} 
on $\R^3\backslash\mathcal{K}$
requires Dirichlet boundary conditions to work.
To get around this,
we note that \eqref{ks2} holds on $\R^3\backslash\mathcal{K}$
provided either $\partial_j$ or the prime is a time derivative,
by using the same argument as in the proof of Lemma \ref{lemmaks}.
This similarly handles the case in which any factor of $Z^\alpha$
involves a time derivative.
For the remaining
cases, involving purely spatial derivatives,
we can use the following
elliptic estimate, which uses the fact that 
$w_k$ vanishes on the boundary, 
\begin{multline*}
\sum_{\substack{|\alpha|\le 9\\ 1\le i,j\le 3}}
\|\partial_x^\alpha \partial_i \partial_j w_{k}(t,\cd)
\|_{L^2(\R^3\backslash \mathcal{K}\,:\,|x|<1)}\\
\leq C\,\sum_{|\alpha|\le 9}\|\partial_x^\alpha \Delta w_{k}(t,\cd)
\|_{L^2(\R^3\backslash \mathcal{K}\,:\,|x|<2)}+
C\,\|w_{k}(t,\cd)\|_{L^6(\R^3\backslash \mathcal{K}\,:\,|x|<2)}\,.
\end{multline*}
%and similarly for $u_0$, since $w_{k-1}$ and $u_0$
%satisfy Dirichlet conditions on $\partial\mathcal{K}$.
The other two terms in \eqref{firstthreeterms} are estimated similarly,
and thus these terms satisfy the bound \eqref{m6}.

Next consider the following terms from the right side of \eqref{m5},
\begin{equation}
\label{43.1}
\sum_{|\alpha|\le 13}
\|\square_\gamma\partial^\alpha_{t,x}w_k(t,\cd)\|_\ltoo
+ \sum_{\substack{|\alpha|+m\le 11\\m\le 1 }}\|\square_\gamma
L^m\partial^\alpha_{t,x}w_k(t,\cd)\|_\ltoo\,.
\end{equation}
We write
\begin{multline*}
(\square_\gamma \partial_{t,x}^\alpha w_k)^I =  \partial_{t,x}^\alpha
\square_\gamma w_k^I -[ \partial_{t,x}^\alpha, \square_\gamma \, ] \, w_k^I \\
= - (1 - \eta) \sum_{\substack{0\le i,j,l\le 3\\ 1\le J,K\le N} }B^{IJ,ij}_{K,l}
\bigl[\partial^\alpha_{t,x}, \partial_l u_0^K+\partial_l w_{k-1}^K \bigr]
\partial_i\partial_j w_k^J\\
+\partial^\alpha_{t,x} (1 - \eta)
\sum_{\substack{0\le i,j,l\le 3\\ 1\le J,K\le N} }B^{IJ,ij}_{K,l}
(\partial_l u_0^K+\partial_l w_{k-1}^K)
\partial_i\partial_j u_0\\
+\partial_{t,x}^\alpha (1 - \eta) B^I\bigl(d(u_0+w_{k-1})\bigr)
-\partial_{t,x}^\alpha [\square_{c_I},\eta]u^I\,.
\end{multline*}
Upon expansion, each term on the right hand side, with the exception of
the last, which has $L^2$ norm bounded by $C\,\e$ by \eqref{10.2}, is the
product of two terms, at least one of which involves at most
8 derivatives. We can obtain $L^\infty$ bounds on such a term
by Sobolev embedding, and hence estimate the $L^2$ norm of the
product by the right hand side of \eqref{m6}. The same argument
applies to the second term in \eqref{43.1}.

Finally, consider
\begin{multline*}
\sum_{|\alpha|\le 12}\|\square_c
\partial^\alpha_{s,x}w_k\|_{L^2([0,T_\e]\times\R^3
\backslash\mathcal{K})}
+\sum_{\substack{|\alpha|+m\le 10\\ m\le 1}}
\|\square_c L^m\partial^\alpha_{s,x}w_k
\|_{L^2([0,T_\e]\times \R^3\backslash\mathcal{K})}
\\
+\sum_{\substack{|\alpha|+m\le 8 \\ m\le 1}}
\|\square_c L^m Z^\alpha w_k\|_{L^2([0,T_\e]\times
\R^3\backslash \mathcal{K})}\,.
\end{multline*}

We write
\begin{multline*}
\square_c\partial_{t,x}^\alpha w_k^I = 
\partial^\alpha_{t,x}
\sum_{\substack{0\le i,j,l\le 3 \\ 1\le J,K\le N}} 
B^{IJ,ij}_{K,l}
(\partial_l u_0^K+\partial_l w_{k-1}^K)
\partial_i\partial_j w_k^J
\\
+\partial^\alpha_{t,x}
\sum_{\substack{0\le i,j,l\le 3 \\ 1\le J,K\le N}}
B^{IJ,ij}_{K,l}
(\partial_l u_0^K+\partial_l w_{k-1}^K)
\partial_i\partial_j u_0
\\
+\partial_{t,x}^\alpha B^I\bigl(d(u_0+w_{k-1})\bigr)
-\partial_{t,x}^\alpha [\square_{c_I},\eta]u^I\,.
\end{multline*}
The last term involving $u$ is easily handled as in \eqref{stuffone} above.
As a representative of the other terms that arise upon expanding
derivatives, noting that $|\alpha|\le 12$, consider
$$
\sum_{|\beta|\le 6,|\mu|\le 13}
\|\partial_{t,x}^\beta w'_{k-1}\,\partial_{t,x}^\mu w'_k
\|_{L^2([0,T_\e]\times
\R^3\backslash \mathcal{K})}\,.
$$
We now apply Lemma \ref{lemma3.3} to the $\beta$ terms and sum over
$R$ to conclude that
\begin{multline*}
\sum_{|\beta|\le 6,|\mu|\le 13}
\|\partial_{t,x}^\beta w'_{k-1}\,\partial_{t,x}^\mu w'_k
\|_{L^2([0,T_\e]\times
\R^3\backslash \mathcal{K})}
\\
\le
C\sum_{|\beta|\le 8}\|\langle x\rangle^{-1}Z^\beta w_{k-1}'
\|_{L^2([0,T_\e]\times
\R^3\backslash \mathcal{K})}
\sup_{0<t<T_\e}\sum_{|\mu|\le 13}\|\partial_{t,x}^\mu w'_k(t,\cd)
\|_\ltoo\\
\le
C\cdot C_1\cdot\e\cdot\ln(T_\e)^{\frac 12}\,M_k(T_\e)\le
C\cdot C_1\cdot\kappa\cdot M_k(T_\e)\,.
\end{multline*}
The other terms are similarly seen to be bounded by the right hand side
of \eqref{m6}, which completes the proof of \eqref{m6}.

Next, using the energy inequality, one observes that  $\{w_k\}$ is
a Cauchy sequence in the energy norm.  Because of this and
\eqref{10.5} we conclude that $w_k$ must converge to a solution of
\eqref{10.3} which satisfies the bounds in \eqref{10.5}.
Consequently, $u=u_0+w$ will be a solution of the original
equation \eqref{0.2} which verifies the analog of \eqref{10.5}.
If the data is $C^\infty$ and satisfies the compatibility
conditions to infinite order, the solution will be $C^\infty$ on
$[0,T_\varepsilon]\times \R^3\backslash \mathcal{K}$ by
standard local existence theory (see e.g., \cite{KSS}).

This completes the proof of Theorem \ref{theorem0.1} .
\qed

\end{document}